\documentclass[preprint,12pt]{elsarticle}

\usepackage{graphicx}
\usepackage{epstopdf, epsfig}
\usepackage{float}

\usepackage{stmaryrd}
\usepackage{amssymb}

\usepackage{pgfplotstable}
\usepackage{amsmath}
\usepackage{graphicx}
\usepackage{natbib}
\usepackage{pgfplots}
\usepackage{tikz}
\usetikzlibrary{shapes}
\usetikzlibrary{arrows}
\usetikzlibrary{arrows.meta}
\usetikzlibrary{patterns}
\usepackage{xcolor}
\usepgfplotslibrary{colormaps}

\usepackage{transparent}

\usepackage{anysize}
\marginsize{2cm}{2cm}{2cm}{2cm}

\usepackage{lineno,hyperref}

\pgfplotsset{compat=newest}
\usepgflibrary{decorations.markings}
\usepgflibrary{decorations.shapes}

\definecolor{green}{rgb}{0.0,.6,0.0}
\definecolor{migris}{rgb}{.85,.85,.85}%

\tikzset{dashed/.style={dash pattern=on 3pt off 1pt}}
\tikzset{dashdot/.style={dash pattern=on .4pt off 3pt on 4pt off 3pt}}
\tikzset{,
	MyPersp1/.style={scale=1.8,x={(-0.8cm,-0.4cm)},y={(0.8cm,-0.4cm)},z={(0cm,1cm)}},
	MyPoints/.style={fill=white,draw=black,thick}
		}
\pgfplotsset{
	axis on top,
	xtick align=center,
	ytick align=center,
	xlabel near ticks,
	ylabel near ticks,
         legend cell align=left,
	width={0.9\textwidth},
    	height={0.0\textwidth},
	scale only axis,
	legend style={font=\tiny},
	title style={font=\footnotesize},
	every axis/.append style={font=\footnotesize},
	ticklabel style={font=\footnotesize},
	xlabel style={font=\footnotesize},
	ylabel style={font=\footnotesize},
	xticklabel style={font=\footnotesize},
	every axis plot/.append style={line width=.75pt,mark size=3pt},
	}

\newcommand{\SigmaR}{\Sigma_R}
\newcommand{\VR}{\mathcal{V}_R}
\newcommand{\V}{\mathcal{V}}
\newcommand{\Ca}{Ca}

\newcommand{\ngrad}{\boldsymbol{n} \! \cdot \! \nabla}

\newcommand{\vzero}{\boldsymbol{0}}
\newcommand{\x}{\times}

\newcommand{\vt}{\boldsymbol{t}}

\newcommand{\vq}{\boldsymbol{q}}
\newcommand{\vn}{\boldsymbol{n}}
\newcommand{\vnS}{\boldsymbol{n}_S}
\newcommand{\ve}{\boldsymbol{e}}

\newcommand{\vv}{\boldsymbol{v}}
\newcommand{\vx}{\boldsymbol{x}}
\newcommand{\vX}{\boldsymbol{X}}

\newcommand{\veps}{\boldsymbol{\varepsilon}}
\newcommand{\vf}{\boldsymbol{f}}

\newcommand{\id}{\mathcal I}
\newcommand{\idS}{\mathcal{I}_S}

\renewcommand{\div}{\nabla \cdot}

\newcommand{\figref}[1]{fig.~\ref{#1}}

\newcommand{\secref}[1]{sec.~\ref{#1}}

\newcommand{\test}[1]{\tilde{#1}}

\newcommand{\stress}{\boldsymbol{\tau}}
\newcommand{\dn}{\rho}
\newcommand{\vdn}{\boldsymbol{\rho}}

\renewcommand{\div}{\boldsymbol{\nabla} \cdot}

\newcommand{\grad}{\boldsymbol{\nabla}}
\newcommand{\hgrad}{\boldsymbol{D}}

\newcommand{\gradS}{\boldsymbol{\nabla}_{\!\!S}}
\newcommand{\hgradS}{\boldsymbol{D}_{\!S}}%\newcommand{\hgradS}{\hat{\nabla}_{\!S}}

\newcommand{\divS}{\boldsymbol{\nabla}_{\!S} \cdot}
\newcommand{\hdivS}{\boldsymbol{D}_{\!S} \cdot}%\newcommand{\hdivS}{\hat{\nabla}_{\!S} \cdot}

\newcommand{\rotS}{\boldsymbol{\nabla}_{\!S}\times}

\newcommand{\dd}{{\rm{d}}}

\newcommand{\rojo}[1]{}

\newcommand{\intV}[2]{ \int_{#2}  #1 \,\dd \V }
\newcommand{\intS}[2]{ \int_{#2}  #1 \,\dd \Sigma }

\newcommand{\vGamma}{\vnS \dd\Gamma}

\usepackage{amssymb}
\journal{ }

\begin{document}

\begin{frontmatter}

%% Title, authors and addresses

%% use the tnoteref command within \title for footnotes;
%% use the tnotetext command for theassociated footnote;
%% use the fnref command within \author or \address for footnotes;
%% use the fntext command for theassociated footnote;
%% use the corref command within \author for corresponding author footnotes;
%% use the cortext command for theassociated footnote;
%% use the ead command for the email address,
%% and the form \ead[url] for the home page:
%% \title{Title\tnoteref{label1}}
%% \tnotetext[label1]{}
%% \author{Name\corref{cor1}\fnref{label2}}
%% \ead{email address}
%% \ead[url]{home page}
%% \fntext[label2]{}
%% \cortext[cor1]{}
%% \address{Address\fnref{label3}}
%% \fntext[label3]{}

\title{PDEs on deformable domains: Boundary Arbitrary Lagrangian-Eulerian (BALE) and Deformable Boundary Perturbation (DBP) methods.}

%% use optional labels to link authors explicitly to addresses:
%% \author[label1,label2]{}
%% \address[label1]{}
%% \address[label2]{}

\author{Javier Rivero-Rodriguez$^{1}$ \& Miguel P\'erez-Saborid$^{2}$ \& Benoit Scheid$^{1}$}

\address{$^{1}$TIPs, Universit\'e Libre de Bruxelles, C.P. 165/67, Avenue F. D. Roosevelt 50, 1050 Bruxelles, Belgium \\
$^{2}$Departamento de Ingenier\'ia Aeroespacial y Mec\'anica de Fluidos, Escuela T\'ecnica Superior de Ingenier\'ia, Universidad de Sevilla, Av. de los Descubrimientos s/n, 41092 Sevilla, Spain}

\begin{abstract}

Many physical problems can be modelled by partial differential equations on unknown domains. Several examples can easily be found in the dynamics of free interfaces in fluid dynamics, solid mechanics or in fluid-solid interactions. To solve these equations in an arbitrary domain with nonlinear deformations, we propose a mathematical approach allowing to track the boundary of the domain, analogue of, and complementary to, the Arbitrary Lagrangian-Eulerian (ALE) method for the interior of the domain. We name this method as the Boundary Arbitrary Lagrangian-Eulerian (BALE) method. Additionally, in many situations nonlinear deformations can be avoided with the help of some analyses which rely on small deformations of the boundary, such as stability analysis, asymptotic expansion and gradient-based shape optimisation. For these cases, we propose an approach to perturb the domain and its boundaries and write the partial differential equations at the unperturbed domain together with the boundary conditions at the unperturbed boundary, instead of at the perturbed ones, which are a priori unknown. We name this method as the Deformable Boundary Perturbation (DBP) method. These two proposed methods rely on the boundary exterior differential operator, whose relevant properties for the present work are evidenced. We show an example for which the BALE and DBP methods are applied, and for which we include the weak formulation revealing the appropriateness of the finite element method in this context. 
 
\end{abstract}

 
\end{frontmatter}

\section{Introduction}

\definecolor{color1}{rgb}{.5,.5,1}
\definecolor{color2}{rgb}{1,.5,.5}

Partial differential equations whose domains are a priori unknown arise in a plethora of physical configurations where geometrical nonlinearities, free interfaces and shape optimisation matter. Geometrical nonlinearities are of crucial importance for the proper description of stability problems involving (i) deformable solids such as beam buckling, first described by Leonhard Euler \cite[]{timoshenko1961theory}, or follower loads \cite[]{patwardhan1999follower}, (ii) solid-gas interactions, such as the flutter of Tacoma Narrows bridge which led to its collapse, the flutter of airwings \cite[]{colera2018numerical}, or cantilevered pipes conveying fluids \cite[]{rivero2015numerical}, and (iii) fluid interfaces with surface tension, such as liquid bridges \cite[]{plateau1873statique}, vibration and break up of liquid domains \cite[] {rayleigh1896theorya} or capillary waves \cite[]{lamb1932hydrodynamics}.

%All this situations takes places in euclidean geometries. However, subdomains are non-euclidean being the most relevant one the boundary of the domain itself. The main characteristic in non-euclidean spaces is the curvature.

There are plenty of examples in the literature in which Cartesian coordinates are used to describe non-Euclidean subdomains \cite[]{paidoussis1998fluid,rubio2013thinnest}, leading to long expressions and tedious algebraic manipulations. Intrinsic coordinates overcome this difficulty and leads to compact expressions that are much easier to manipulate. These coordinates also facilitate the description of differential operators as we have shown in \cite[]{rivero2015numerical,rivero2017efficient}. However, intrinsic coordinates require to be related to a reference system, thus appealing for a method that describes intrinsic operators in that reference frame.

Physical problems in deformable geometries are usually modelled for simple ones, which can be analytically described in common coordinate systems, such as cylindrical \cite[]{tomotika1935instability}, spherical \cite[]{rayleigh1896theorya} or toroidal \cite[]{penas2015dissolution}, or modifications of the latter such as cylindrical with scaling of the radial direction depending on the axial position \cite[]{herrada2016numerical}. However, such coordinate system are not always easy, practical nor even possible to find for complex geometries, and a change of variables that fulfils a given partial differential equation is required. This is known as arbitrary Lagrangian-Eulerian (ALE) method, originally developed by \citep{Noh64,Fra64}, and has been widely used in the literature for fluid-structure interaction problems \citep{donea1982arbitrary,zhang2001analysis} and for fluid-mechanical problems with moving interfaces \cite[]{ramaswamy1987arbitrary,basaran1992nonlinear,chen1993nonlinear,souli2001arbitrary,balestra2018viscous}. Although a simple Lagrangian description is usually enough for solid mechanics problems with large deformation, ALE has also been used in this field, such as in the crack propagation \cite{koh1988dynamic} or the metal forming processes \cite{gadala2002mesh}. A survey of the method has been recently provided by its author and coworkers \cite{donea2017arbitrary}. However, ALE requires boundary conditions and thus translates the difficulty to the boundary. Although many analytical alternatives are used in the literature, such as following the boundary in a Lagrangian framework \cite[]{basaran1992nonlinear}, in the normal direction to the boundary \cite[]{balestra2018viscous} or in one of the directions of the reference system \cite[]{herrada2016numerical}, there is a lack of a systematic method to avoid the associated drawbacks of remeshing. Following the same spirit as ALE, we propose in this work a novel systematic method to obtain the change of variables needed at the boundary. We named this method as Boundary Arbitrary Lagrangian Eulerian (BALE) method.
 
Analytical changes of variables and description at the boundary further allow to carry out linearisation of model equations in canonical geometries \cite{rayleigh1896theorya,lamb1932hydrodynamics,chandrasekhar1961hydrodynamic}, and in more involved geometries such as in stretched axisymmetric jets \cite[]{herrada2016numerical,rubio2013thinnest}, curve jets and pipes \cite[]{arne2018whipping,rivero2015numerical}, planar plates \cite[]{colera2018numerical} or coflowing liquids \cite[]{castro2012slender}, with the help of methods developed on purpose for those specific geometries. Despite this concept was developed shortly after the development of perturbation techniques \cite{feshbach1941perturbation}, they are either analytical \cite{feshbach1941perturbation,dancer1997domain} or they depend on the numerical method used to solve the equations. In addition, the perturbation leads to an increase of the order of the highest derivative present in the equation, which requires the use of higher order test functions in the finite element method \cite{givoli2000boundary}, or of hypersingular kernels in the boundary integral method \cite{minutolo2004shape}. 

In order to avoid the previous drawbacks, we have devised a systematic method to carry out the linearisation of the equations corresponding to different mathematical models formulated in arbitrary domains, which is independent of the numerical resolution method. Our method, which we have named Deformable Boundary Perturbation (DBP), is based on partial differential equations written for the unperturbed domain and its boundaries, instead of using the a priori unknown perturbed or linearised ones. 

As we shall see, the BALE and DBP methods are based on the boundary exterior differential operator, for which we have provided a brief and self-consistent description in \secref{Exdef}, where we establish its relation with nabla operators and the Stokes theorem. In \secref{BALE}, we introduce the BALE method. In \secref{DBP}, we introduce the DBP method and apply it to the perturbation of integrals in \secref{PInt}, to the treatment of mixed boundary conditions in \secref{PMix} and to the boundary exterior differential operator in \secref{AppD}. In \secref{Linearisation}, we show how to linearise a system of PDE defined on deformable domain first using the DBP method and then introducing a regular expansion. In \secref{Example1}, we provide with an example of a system of partial differential equations on a deformable domain in \secref{zeroth}, its geometry perturbation and linearisation in \secref{lin}, the weak formulation appropriate for the finite element method in \secref{WF} and the discussion of the results provided by the BALE and DBP methods in \secref{Disc}. Finally, in \secref{Conclusions}, we present the conclusions.

\section{Differential operators}   \label{Exdef}  
     
In this section, we introduce the differential operators in Euclidean spaces, which are necessary for this work, namely the directional derivative, the nabla operator and the exterior differential operator, as well as their boundary counterparts. We also introduce the boundary Stokes theorem and the boundary reciprocal theorem. But one first needs to introduce the set or subset of the space at which the differential operators are defined.

We consider a domain, denoted $\V$, with identity tensor $\id$, embedded in an Euclidean $3$ space. The boundary of the domain, denoted by $\Sigma \equiv \partial \V$ $\Sigma$, is contained in the same Euclidean space as $\V$ and its outer normal will be denoted as $\vn$. The boundaries, with identity tensor $\idS = \id - \vn \vn$, can be also considered as embedded in an Euclidean $2$ sub-space. The contour of a boundary, denoted by $\Gamma \equiv \partial  \Sigma$, is contained in the same Euclidean sub-space as $\Sigma$ and its boundary outer normal will be denoted as $\vnS$, such as $\vn \cdot  \vnS = 0$ and $\idS \cdot \vnS = \vnS$. We use the laboratory frame $\{ \ve_i \}_{i=1,2,3}$, in which the position vector writes $\vx=x_i \ve_i$ with coordinates $x_i$. Index $i$ gets values $1$, $2$ and $3$, as well as the latin indices $j$ and $k$, whereas greek indices $\alpha$ and $\beta$ get values $1$ and $2$.
 
 \subsection{Definitions}
 
The directional derivative $\partial_{x_i}$ and boundary directional derivative $\partial_{S x_i}$ are defined as
\begin{subequations}\label{dirder}
\begin{align}
\partial_{x_i} \varphi (\vx) &= \partial_{\epsilon} \varphi(\vx+\epsilon \ve_i ) \big \vert_{\epsilon=0} \,, \\
\partial_{S x_i} \varphi (\vx) &= \partial_{\epsilon} \varphi(\vx+\epsilon \idS \cdot \ve_{i} ) \big \vert_{\epsilon=0} \,, \label{dirderSurf}
\end{align}
\end{subequations}
where $\varphi$ is a generic quantity which can be a scalar, a vector or a tensor. Using the chain rule, both operators in \eqref{dirder} can be related by
\begin{align}
\label{Sii}
\partial_{S x_i} \varphi = \ve_i \cdot \idS \cdot \ve_j \partial_{x_j} \varphi \,,
\end{align}
 as obtained from substituting $\varphi= \ve_j \cdot \vx$ in \eqref{dirderSurf}, and with the use of Einstein's convention.

At the implementation level of numerical schemes to solve partial differential equations, only directional derivatives appear. For this reason, the differential operators will be defined in terms of directional derivatives. The nabla operator $\grad$ and boundary nabla operator $\gradS$ are defined as
\begin{subequations}
\label{nablaop}
\begin{align} 
\label{grad}
\grad \varphi &=  \ve_i \partial_{x_i} \left(  \varphi \right)  \,, \\
\label{gradS}
\gradS \varphi &= \ve_i \partial_{S x_i} \left(    \varphi \right)  \,.
\end{align}
\end{subequations}
where $\gradS$ can also be written, after substituting \eqref{Sii} in \eqref{gradS}, as $\gradS = \idS \cdot \grad$.
Thus, according to \eqref{Sii} and \eqref{gradS}, the boundary partial derivative $\partial_{Sx_i} \varphi $ is just the $\ve_i$ component of $\gradS \varphi$ and, the operator $\gradS$ can be geometrically interpreted as the projection of the operator $\grad$ onto the boundary.  

Likewise, the exterior differential operator $\hgrad$ and boundary exterior differential operator $\hgradS$ are defined as
\begin{subequations}
\label{hgrad_hgradS}
\begin{align}
\hgrad   \varphi &= \ve_i \cdot \partial_{x_i} \left( \id  \varphi \right)  \,, 
\label{hgrad} \\
\hgradS \varphi &= \ve_i \cdot \partial_{S x_i} \left(  \idS  \varphi \right) \,,
\label{hgradS}
\end{align}
\end{subequations}
in which, as compared to the nabla operators, the laboratory frame is projected onto the identity tensor of the space or sub-space where the derivative is defined. Substituting \eqref{nablaop} in \eqref{hgrad_hgradS} leads to
\begin{subequations}
\label{refertonabla}
\begin{align}
%%%\gradS   \varphi &= \left(\idS \cdot \grad\right) \varphi \label{gradS_geo}  \,, \\
\hgrad   \varphi &= \grad   \varphi  \,, \\
\hgradS \varphi &=  \gradS \cdot \left( \idS \varphi  \right) \,,
\end{align}
\end{subequations}
where the nabla operator and the exterior differential operator are equivalent since $\id \cdot \ve_i = \ve_i  $. Hence, no distinction will be made hereafter, and either $\hgrad$ or $\grad$ will be used for the sake of analogy to their boundaries counterparts. However, since $\idS \cdot \ve_i \neq \ve_i  $ in general, boundary operators are not equivalent and are related by
\begin{align}
\label{equiv}
\hgradS   \varphi  = \gradS \varphi + \left( \gradS \cdot \idS \right)  \varphi  \,,
\end{align}
where $\gradS \cdot \idS = - \vn \gradS \cdot \vn$ is the mean curvature vector of the boundary.

\subsection{Stokes theorem}

The geometrical interpretations of $\hgrad $, or equivalently $\grad$, and $\hgradS$ come out of the Stokes theorems. Starting with the $\hgrad$ operator, which can be written in any curvilinear coordinates $\xi_i$ such as $\vx=\vx(\xi_i)$, and applying the chain rule $\partial_{ x_i} = \left( \partial_{ x_i} \xi_{j} \right) \partial_{\xi_{j}}$ to \eqref{hgrad}, we obtain
\begin{equation}
\label{eq7}
\hgrad \varphi = \grad \xi_i \partial_{\xi_i} \varphi\,. 
\end{equation}
The Stokes theorem comes out by integrating \eqref{eq7} over a domain $ \V $, such as $ \xi_i \in [\xi_{i}^{-}, \xi_{i}^{+}]$, and bounded by $\Sigma$, which can be decomposed into the boundaries $\Sigma_{i}^{\pm}$ at coordinates surfaces $\xi_{i} = \xi_{i}^{\pm}$. In effect, this integral writes as
\begin{align}
\label{eq8}
\int_{\V} \hgrad \varphi \, \dd \V  =
\int_{\xi_{1}^{-}}^{\xi_{1}^{+}} \int_{\xi_{2}^{-}}^{\xi_{2}^{+}} \int_{\xi_{3}^{-}}^{\xi_{3}^{+}} 
J  \, \grad \xi_i  \,  \partial_{\xi_i} \varphi
  \, \dd \xi_{3}\, \dd \xi_{2}\, \dd \xi_{1} 
 \,  ,
\end{align}
where $J= \left( \partial_{\xi_1} \vx \x \partial_{\xi_2} \vx \right) \cdot \partial_{\xi_3} \vx $ is the Jacobian.
Taking into account that
\begin{subequations}
\begin{align}  
J \, \grad \xi_i   &=\frac12 \epsilon_{ijk} \partial_{\xi_j} \vx \x \partial_{\xi_k} \vx \,, \\
\partial_{\xi_i} ( J \, \grad \xi_i & ) = \vzero \,,
\end{align}
\end{subequations}
since $\grad \vx = \id = \partial_{\xi_i} \vx \grad \xi_i$, and where $\epsilon_{ijk}$ is the Levi-Civita symbol, \eqref{eq8} can be rewritten as
\begin{align}
\label{VolGenStokes0}
\int_{\V} \hgrad \varphi \, \dd \V  =
\int_{\xi_{1}^{-}}^{\xi_{1}^{+}} \int_{\xi_{2}^{-}}^{\xi_{2}^{+}} \int_{\xi_{3}^{-}}^{\xi_{3}^{+}} \partial_{\xi_i} \Big{[}  \frac12 \epsilon_{ijk} \left( \partial_{\xi_j} \vx \x \partial_{\xi_k} \vx \right) \varphi  \Big{]}
  \, \dd \xi_{1}\, \dd \xi_{2}\, \dd \xi_{3} 
 \,  ,
\end{align}
which, using the fundamental theorem of calculus and identifying the differential boundary outer normal vector $\vn d \Sigma$ of the boundary $\Sigma_{i}^{\pm}$ as $ \pm \frac12 \epsilon_{ijk} \left( \partial_{\xi_j} \vx \x \partial_{\xi_k} \vx \right) \dd\xi_j \dd \xi_k$, leads to the Stokes theorem
\begin{align}
\label{VolGenStokes}
\int_{\V} \hgrad \varphi \, \dd \V  =
\int_{\Sigma} \vn \varphi \, \dd \Sigma
 \,  .
\end{align}

%%%%%%%%
In the same spirit, the $\hgradS$ operator can be written in any boundary curvilinear coordinates $\xi_{\alpha}$ such as $\vx=\vx(\xi_\alpha)$, and applying the chain rule $\partial_{S x_i} = \left( \partial_{S x_i} \xi_{\alpha} \right) \partial_{\xi_{\alpha}}$ to \eqref{hgradS}, we obtain
\begin{equation}
\label{eq12}
\hgradS \varphi = \gradS \xi_\alpha \cdot \partial_{\xi_\alpha} \left( \idS \varphi \right) \,.
\end{equation}
The boundary Stokes theorem comes out by integrating \eqref{eq12} over a boundary $ \Sigma $, such as $ \xi_\alpha \in [\xi_{\alpha}^{-}, \xi_{\alpha}^{+}]$, and contoured by $\Gamma$, which can be decomposed into the contours $\Gamma_{\alpha}^{\pm}$ at coordinates lines $\xi_{\alpha} = \xi_{\alpha}^{\pm}$. In effect, this integral writes as
\begin{align}
\label{eq13}
\int_{\Sigma} \hgradS \varphi \, \dd \Sigma  =
\int_{\xi_{1}^{-}}^{\xi_{1}^{+}} \int_{\xi_{2}^{-}}^{\xi_{2}^{+}}  
J_S   \, \gradS \xi_\alpha  \cdot \partial_{\xi_\alpha} \left( \idS \varphi \right)
  \, \dd \xi_{2}\, \dd \xi_{1}
 \,  ,
\end{align}
where $J_S = \left( \partial_{\xi_1} \vx \x \partial_{\xi_2} \vx \right) \cdot \vn $ is the boundary Jacobian. Taking into account that
\begin{subequations}
\begin{align}
  J_S \, \gradS \xi_\alpha  &=\frac12 \epsilon_{\alpha \beta 3} \partial_{\xi_\beta} \vx \x \vn \,,\\
\idS \cdot \partial_{\xi_\alpha} (   J_S \, \gradS \xi_\alpha & ) = \vzero \,,
\end{align}
\end{subequations}
since $\gradS \vx = \idS = \partial_{\xi_\alpha} \vx \gradS \xi_\alpha$, and vectors $\partial_{\xi_{\alpha}} \vn$ and $\partial_{\xi_{\beta}} \vx$ are contained in the boundary, \eqref{eq13} can be rewritten as
\begin{align}
\label{SurfGenStokes0}
\int_{\Sigma} \hgradS \varphi \, \dd \Sigma  =
\int_{\xi_{1}^{-}}^{\xi_{1}^{+}} \int_{\xi_{2}^{-}}^{\xi_{2}^{+}} \partial_{\xi_\alpha} \Big{[}  \frac12 \epsilon_{\alpha \beta 3} \left( \partial_{\xi_\beta} \vx \x \vn \right) \varphi  \Big{]}
  \, \dd \xi_{2}\, \dd \xi_{1}
 \,  ,
\end{align}
which using the fundamental theorem of calculus and identifying the differential contour outer normal vector $\vnS d \Gamma$ of the contour $\Gamma_{\alpha}^{\pm}$ as $ \pm \frac12 \epsilon_{\alpha\beta3} \left( \partial_{\xi_\beta} \vx \x \vn \right) \dd\xi_\beta$, leads to the boundary Stokes theorem
\begin{align}
\label{SurfGenStokes}
\int_{\Sigma} \hgradS \varphi \, \dd \Sigma  =
\int_{\Gamma} \vnS \varphi \, \dd \Gamma
 \,  .
\end{align}
The trace of \eqref{SurfGenStokes} for a vector or tensor quantity $\varphi$ is known as the surface divergence theorem \cite[page 239]{weatherburn2016differential}.

Note that the theorems \eqref{VolGenStokes} and \eqref{SurfGenStokes} are independent of the chosen curvilinear system whose use is only convenient for demonstration purposes.

\subsection{Reciprocal theorem}

The discretisation of \eqref{refertonabla} using the finite element method can be facilitated with the use of the reciprocal theorems. Multiplying \eqref{refertonabla} by an arbitrary scalar test function $\psi=\psi(\vx)$ and after application of the chain rule using \eqref{nablaop} and \eqref{hgrad_hgradS}, one obtains the reciprocal theorem and boundary reciprocal theorem in differential form,
\begin{subequations}
\label{recip_diff}
\begin{align}
\psi \hgrad \varphi  &= \hgrad   ( \psi \varphi) -  (\grad \psi) \varphi   \,, \\
\psi \hgradS \varphi  &= \hgradS   ( \psi \varphi) -  (\gradS \psi) \varphi   \,.
\end{align}
\end{subequations}
Integrating \eqref{recip_diff} over an arbitrary domain $\V$ or boundary $\Sigma$, respectively, leads to the reciprocal theorem and boundary reciprocal theorem in integral form
\begin{subequations}
\label{prop3}
\begin{align}
\int_\V \psi \hgrad \varphi \dd \V &= \int_\Sigma \vn \psi \varphi  \dd \Sigma - \int_\V  (\grad \psi)  \varphi \dd \V  \,,\\
\int_\Sigma \psi \hgradS \varphi \dd \Sigma &= \int_\Gamma \vnS \psi \varphi  \dd \Gamma - \int_\Sigma  (\gradS \psi)  \varphi \dd \Sigma  \,,  \label{prop3B}
\end{align}
\end{subequations}
where the Stokes theorems \eqref{VolGenStokes} and \eqref{SurfGenStokes} have been used. It should be noted that \eqref{prop3B} provides a useful manner to discretise $\hgradS$, which avoids the derivative of $\idS$ appearing in \eqref{hgradS}.

%%%%%%%%%%%%%%%%%%
%%%%%%%%%%%%%%%%%%
%%%%%%%%%%%%%%%%%%
%%%%%%%%%%%%%%%%%%

\section{Boundary Arbitrary Lagrangian-Eulerian (BALE) method}  \label{BALE} 

\begin{figure}[h!]
\centering
\input{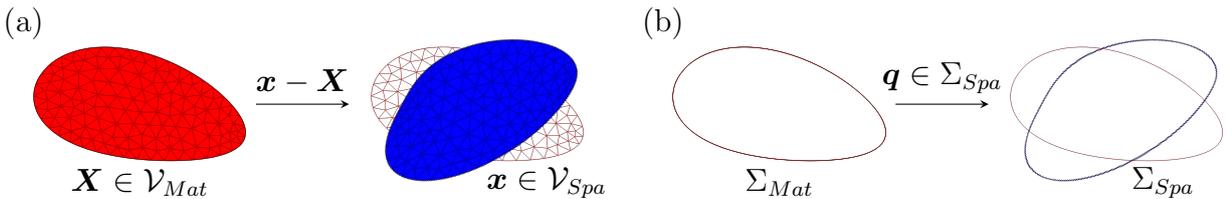}
\caption{Scheme of (a) domain and (b) boundary deformation, which is used in the BALE method.}
\label{SketchBALEvsALE}
\end{figure}

In this section, we describe a procedure to model the nonlinear deformation of a domain and its boundary, as illustrated in \figref{SketchBALEvsALE}, which consists in the ALE method to deform the domain and the BALE method to deform the boundary. The undeformed domain is usually referred as material, $\V_{Mat}$, and the deformed domain as spatial, $\V_{Spa}$. The position vectors within the domains are $\vX\in \V_{Mat}$ and  $\vx\in \V_{Spa}$, respectively, such as the material domain $\V_{Mat}$ maps into the spatial domain $\V_{Spa}$ under the effect of the displacement field $\vx-\vX$. Since domains are exclusively defined by their boundaries, $\Sigma_{Mat} \equiv \partial \V_{Mat} $ and $\Sigma_{Spa} \equiv \partial \V_{Spa} $, and an orientation of the outer normal vector, any displacement field $\vx-\vX$ with the same displacement of the boundary, $\vq$, i.e. 
 \begin{equation}
 \label{Sigma0}
\vx - \vX = \vq  \qquad \mbox{at } \Sigma_{Spa} \,,
\end{equation}
leads to the same domain. It is worth mentioning that the differential operators introduced in \secref{Exdef} are defined on the spatial domain and spatial boundary, i.e. on $\vx$, although they can also be defined on the material ones after the change of variables $\vx-\vX$. The material domain is discretised and leads to different discretisation of the spatial domain, depending on the choice of $\vx-\vX$ 

There is a infinite number of $\vx - \vX$ and we list here a few. The first one is to follow the material point where $\vx - \vX$ is set to the displacement of the material point. This method is usually not appropriate for boundary conditions which are naturally described in an Eulerian framework and another choice becomes mandatory \cite{huerta1988viscous}. An analytical displacement field $\vx - \vX$ can be chosen such as it follows the boundary of the domain as done by \cite{herrada2016numerical}. However, it is not a general solution since an analytical transformation is not always easy, possible nor practical to find, and thus, its use is limited to simple geometries. To overcome this difficulty, one can choose a transformation that fulfils a given partial differential equation (PDE), which can be conveniently chosen to minimise the loss of mesh quality. The boundary conditions are then used to adequately follow the boundaries. The latter method is known as the Arbitrary Lagrangian-Eulerian (ALE) method \cite{hirt1974arbitrary}. We choose the PDE to be the Laplace equation, although other choices are possible such as elasticity equation \cite{souli2001arbitrary},
 \begin{equation}
 \label{ALE}
\div \grad (\vx - \vX) = \vzero \qquad \mbox{at } \V_{Spa} \,,
\end{equation}
together with the boundary condition \eqref{Sigma0},

If the displacement of the boundary is within itself, i.e. spatial and material boundaries coincide, spatial and material domains also coincide. Since the boundary is a $2$ sub-space, it reduces the degrees of freedom of $\vq$ by $2$, i.e. to $1$, which requires the choice of the remaining one. Several choices are possible similar to the case of the domain. Like for the transformation inside the domain, one can also use the methods described in the previous paragraph. The first one is to follow the material point, which has some drawbacks such as large deformation which usually requires remeshing \cite{balestra2018viscous} as well as the impossibility of carrying out stationary analysis if Eulerian boundary conditions are imposed at deformable interfaces as it occurs in fluid interfaces subject to surface tension \cite{chen2014inertial,Rivero2018a}. Then, to avoid the latter analytical displacement can be an appropriate alternative as used by \cite{herrada2016numerical} who choose the displacement to be in the direction of one vector of the coordinate system. Other choices such as displacement in the normal direction are also possible. However, these solutions are not general and are limited to simple geometries. To overcome this difficulty, one can choose a transformation of the boundary that fulfils a given boundary partial differential equation (BPDE) with an additional unknown variable which represents the degree of freedom of the boundary. Again, the BPDE can be chosen to minimise the loss of mesh quality, similarly to the case of ALE. In this work, we choose the boundary transformation that fulfill the boundary Poisson equations
\begin{equation}\label{BALEeq}
\hdivS \gradS \vq = g \,\vn \qquad \mbox{at } \Sigma_{Spa} \,,
\end{equation}
where the source terms has only one component with magnitud $g$ in the direction perpendicular to the boundary which represents the degree of freedom of the boundary. It is worth noting that \eqref{BALEeq} is the Poisson equation out of the boundary $\Sigma_{Spa}$ and Laplace equation within it, being the latter analogous to the Laplace equation for the domain \eqref{ALE}. Despite we have considered that $\vq$ is a finite displacement, it can also be used for infinitesimal displacements.

%%This choice is not unique and other BPDEs such as
%%\begin{equation}
%%\hdivS \gradS \vq_S = \hgradS g  \qquad \mbox{at } \Sigma \,,
%%\end{equation}
%%also lead to satisfactory results.

\section{Deformable Boundary Perturbation (DBP) method}  \label{DBP}  
          
\begin{figure}[h!]
\centering
\input{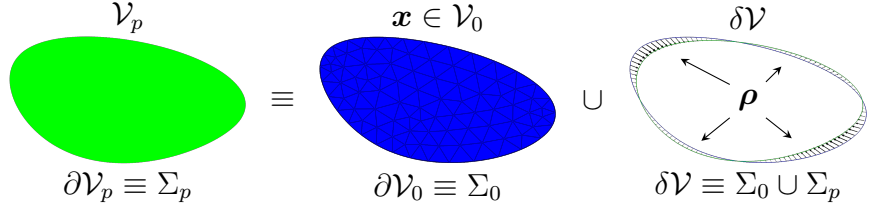}
\caption{Scheme of a planar section of the decomposition of a perturbed domain and its boundary, which is used in the DBP method.}
\label{DomDecomp}
\end{figure}

Let us consider a perturbed domain $\V_p$ with perturbed boundary $\Sigma_p$ as schematised in \figref{DomDecomp}, which comes out of an infinitesimal displacement $\vdn$ of the unperturbed boundary $\Sigma_0$. Sweeping the unperturbed boundary $\Sigma_0$ along an infinitesimal displacement $\vdn$ generates the perturbation domain $\delta V$, where $\delta$ denotes variation between perturbed and unperturbed domain. As shown in \figref{DomDecomp}, the perturbed domain $\V_p$ can be decomposed in the unperturbed domain $\V_0$ and the perturbation domain $\delta \V$. The goal of the DBP method is to rewrite in the unperturbed domain and boundary, the equations initially defined on the perturbed ones. For this reason, only the unperturbed domain and boundary need to be discretised.
 
\begin{figure}[h!]
\centering
\begin{tikzpicture}[scale=5]

\node at (0.2,1.165) {(a)};

\draw [blue,fill=blue!50!white,thick] plot [smooth cycle,tension=0.5] coordinates {
(0.75622,0.51496) 
(0.67788,0.59073) 
(0.55115,0.64902) 
(0.44286,0.65902) 
(0.32765,0.63445) 
(0.22553,0.57325) 
(0.19862,0.50913)
(0.21475,0.42461)
(0.30000,0.36758)
(0.36452,0.38089)
(0.42903,0.43335)
(0.48664,0.46833)
(0.55115,0.45958)
(0.61567,0.41587)
(0.66406,0.38964)
(0.72627,0.38964)
(0.77005,0.41587)
(0.78018,0.45084)
};

\draw [green,fill=green!50!white,thick] plot [xshift=1.5,yshift=9.5,rotate=10,scale=1,smooth cycle,tension=0.5] coordinates {
(0.75622,0.51496) 
(0.67788,0.59073) 
(0.55115,0.64902) 
(0.44286,0.65902) 
(0.32765,0.63445) 
(0.22553,0.57325) 
(0.19862,0.50913)
(0.21475,0.42461)
(0.30000,0.36758)
(0.36452,0.38089)
(0.42903,0.43335)
(0.48664,0.46833)
(0.55115,0.45958)
(0.61567,0.41587)
(0.66406,0.38964)
(0.72627,0.38964)
(0.77005,0.41587)
(0.78018,0.45084)
};

\draw [orange,thick,->] plot  coordinates {
(0.66406,0.38964)
(0.64406,0.83)
};

\draw [orange,thick,->] plot  coordinates {
(0.77005,0.41587)
(0.74005,0.874)
};

\fill[opacity=.5,fill=orange!50!white,line width=2]  
plot coordinates{(0.66406,0.38964)(0.72627,0.38964)(0.77005,0.41587) }
--	
plot coordinates{(0.74005,0.874)(0.69005,0.835)(0.64406,0.83)}
-- cycle
;

\node[green] at (0.2,0.7) {$C_p$};
\node[green] at (0.4,0.95) {$S_p$};
\node[blue] at (0.5,0.7) {$C_0$};
\node[blue] at (0.4,0.55) {$S_0$};

\node[orange] at (0.84,0.72) {$\dd S_g$};

\node[orange] at (0.7,0.32) {$\dd C_0$};

\end{tikzpicture}
\begin{tikzpicture}[x=.42\textwidth,y=.42\textwidth,scale=0.92]
	%-{Latex[length=2mm]}
	%\node at (-1,1.5)   {$(b)$};
	
	\node at (-0.8,1.5) {(b)};
	
	\draw[dashed]  
	plot [domain = 60:120, samples = 40, variable = \i]  ({(1+.0*\i/100)*cos \i}, {(1+.0*\i/100)*sin \i}, 0); 
	\draw  
	plot [domain = 110:60, samples = 40, variable = \i]  ({(1+.4*\i/100)*cos \i}, {(1+.4*\i/100)*sin \i}, 0) ;

%	Relleno
	\draw[black!20!white,fill=black!20!white,line width=2]  
	plot [domain = 75:105, samples = 40, variable = \i]  ({(1+.0*\i/100)*cos \i}, {(1+.0*\i/100)*sin \i}, 0) 
	--	
	plot [domain = 100:70, samples = 40, variable = \i]  ({(1+.4*\i/100)*cos \i}, {(1+.4*\i/100)*sin \i}, 0) 
	-- cycle
	;

	\draw[blue!50!white,line width=3]  
	plot [domain = 75:105, samples = 40, variable = \i]  ({(1+.0*\i/100)*cos \i}, {(1+.0*\i/100)*sin \i}, 0) ;
	\draw[green!50!white,line width=3]  
	plot [domain = 100:70, samples = 40, variable = \i]  ({(1+.4*\i/100)*cos \i}, {(1+.4*\i/100)*sin \i}, 0) ;

	\draw[orange!50!white,line width=3] ({cos 75},{sin 75}) -- ({(1+.4*.7)*cos 70},{(1+.4*.7)*sin 70});
	\draw[orange!50!white,line width=3] ({cos 105},{sin 105}) -- ({(1+.4*1.0)*cos 100},{(1+.4*1.0)*sin 100});

	\begin{scope}[shift={(cos 75, sin 75)}]
		\draw[->,blue,thick] (0,0)--({.1*cos 15}, {.1*sin(-15)}) node[below] {$\vn_{S0}$};
		\draw[->,blue,thick] (0,0)--({.1*sin 15}, {.1*cos 15}) node[left] {$\vn_0$};
	\end{scope}
	\begin{scope}[shift={(cos 105, sin 105)}]
		\draw[->,blue,thick] (0,0)--({-.1*cos 15}, {-.1*sin 15 }) node[below] {$\vn_{S0}$};
		\draw[->,blue,thick] (0,0)--({.1*sin -15}, {.1*cos 15}) node[left] {$\vn_0$};
	\end{scope}

%%%%
	\begin{scope}[shift={({(1+.4*0.7)*cos(70)},{(1+.4*0.7)*sin(70)})}]
		\draw[->,green,thick] (0,0)--({.1*cos 33}, {.1*sin(-33)}) node[below] {$\vn_{Sp}$};
		\draw[->,green,thick] (0,0)--({.1*sin 33}, {.1*cos 33}) node[right] {$\vn_p$};
	\end{scope}
	\begin{scope}[shift={({(1+.4*1.0)*cos(100)},{(1+.4*1.00)*sin(100)})}]
		\draw[->,green,thick] (0,0)--({-.1*cos 2}, {-.1*sin(2)}) node[above] {$\vn_{Sp}$};
		\draw[->,green,thick] (0,0)--({.1*sin -2}, {.1*cos 2}) node[right] {$\vn_p$};
	\end{scope}

	\begin{scope}[shift={({.5*cos(105)+.5*(1+.4*1.0)*cos(100)},{.5*sin(105)+.5*(1+.4*1.0)*sin(100)})}]
	
	%	\draw[->,red,thick] (0,0)--({.2*cos(15)*cos(22)}, {.2*sin(-15)*cos(22)}) node[right] {$\vn_{0}\vn_{S0} \cdot \vdn \dd \Gamma$};
	%	\draw[->,red,thick] (0,0)--({-.2*sin(15)*sin(22)}, {-.2*cos(15)*sin(22)}) node[right] {$\vn_{0}\vn_{S0} \cdot \vdn \dd \Gamma$} ;
	
		\draw[->,orange,thick] (0,0)--({-.3*cos(15)*cos(17)}, {-.3*sin(15)*cos(17)}) node[below] {$ \vdn \cdot \vn_{0} \vn_{S0}$};
		\draw[->,orange,thick] (0,0)--({.3*sin(-15)*sin(17)}, {.3*cos(15)*sin(17)}) node[anchor=south east] {$- \vdn \cdot \vn_{S0} \vn_0 $};
		
		\begin{scope}[shift={({.3*sin(-15)*sin(17)}, {.3*cos(15)*sin(17)})}]
		\draw[orange,dashed] (0,0)--({-.3*cos(15)*cos(17)}, {-.3*sin(15)*cos(17)});
		\end{scope}
		\begin{scope}[shift={({-.3*cos(15)*cos(17)}, {-.3*sin(15)*cos(17)})}]
		\draw[orange,dashed] (0,0)--({.3*sin(-15)*sin(17)}, {.3*cos(15)*sin(17)});
		\end{scope}
	
	\begin{scope}[rotate=-2]
		\draw[->,very thick,orange] (0,0)--(-.3,0) node[anchor=south,pos=1.2] {$\vdn \!\!\x\!\! (\vn_{S0} \!\x\! \vn_0) \dd C_0$};
	\end{scope}
		
	\end{scope}

	\begin{scope}[shift={({.5*cos(75)+.5*(1+.4*.7)*cos(70)},{.5*sin(75)+.5*(1+.4*.7)*sin(70)})}]
	\begin{scope}[rotate=-37]
		\draw[->,very thick,orange] (0,0)--(.2,0) node[right] {$\vdn \!\x\! (\vn_{S0} \!\x\! \vn_0) \dd C_0$};
	\end{scope}
%		\draw[->,red,thick] (0,0)--({.2*cos(15)*cos(22)}, {.2*sin(-15)*cos(22)}) node[right] {$\vn_{0}\vn_{S0} \cdot \vdn \dd \Gamma$};
%		\draw[->,red,thick] (0,0)--({-.2*sin(15)*sin(22)}, {-.2*cos(15)*sin(22)}) node[right] {$\vn_{0}\vn_{S0} \cdot \vdn \dd \Gamma$} ;
	
	%	\draw[->,red,thick] (0,0)--({.1*cos 15}, {.1*sin(-15)}) node[above] {$\vn_{S0}$};
	%	\draw[->,red,thick] (0,0)--({.1*sin 15}, {.1*cos 15}) node[left] {$\vn_0$};
	\end{scope}

	\draw[-{Latex[length=2.5mm]},thick] ({cos 75},{sin 75}) -- ({(1+.4*.7)*cos 70},{(1+.4*.7)*sin 70});
	\draw[-{Latex[length=2.5mm]},thick] ({cos 82.5},{sin 82.5}) -- ({(1+.4*.775)*cos 77.5},{(1+.4*.775)*sin 77.5});
	\draw[-{Latex[length=2.5mm]},thick] ({cos 90},{sin 90}) -- ({(1+.4*.85)*cos 85},{(1+.4*.85)*sin 85}) node[midway,left] {$\vdn$};
	\draw[-{Latex[length=2.5mm]},thick] ({cos 97.5},{sin 97.5}) -- ({(1+.4*.925)*cos 92.5},{(1+.4*.925)*sin 92.5});
	\draw[-{Latex[length=2.5mm]},thick] ({cos 105},{sin 105}) -- ({(1+.4*1.0)*cos 100},{(1+.4*1.0)*sin 100});

	\node at ({cos 75},{sin 75}) [circle,fill=blue,inner sep=2pt]{};
	\node at ({(1+.4*.7)*cos 70},{(1+.4*.7)*sin 70}) [circle,fill=green,inner sep=2pt]{};
	\node at ({cos 105},{sin 105}) [circle,fill=blue,inner sep=2pt]{};
	\node at ({(1+.4*1.0)*cos 100},{(1+.4*1.0)*sin 100}) [circle,fill=green,inner sep=2pt]{};

%%%        \node at (.3,1.4)   {$\dd \Sigma_{0} \times \dn $};
%%%    \node at (.3,1.35)   {$\dd \Sigma_{0} $};
%%%     \node at (0,0.9)   {$\dd \Sigma_{0} $};
        
        \begin{scope}[shift={(0,1)}]
	\draw[->,blue, very thick] (0,0)--(0,-.1) node[below] {$-\vn_0 \dd \Sigma_0$};
        \end{scope}
        
        \begin{scope}[shift={ ({(1+.4*85/100)*cos 85}, {(1+.4*85/100)*sin 85}, 0) }]
	\begin{scope}[rotate=-14.5]
	\draw[->,green, very thick] (0,0)--(0,.120) node[right] {$\vn_{p} \dd \Sigma_p$};
        \end{scope}
        \end{scope}
        
        \node at (-.15,1.26)   {$U$};
        
        \node[green] at (-.05,1.41)   {$S_p$};
        \node at (-.10,0.95)   {$S_0$};
        
        \node[green] at (-.18,1.41)   {$C_p$};
        \node at (-.24,0.92)   {$C_0$};
        \node[green] at (+.42,1.27)   {$C_p$};
        \node at (+.25,0.92)   {$C_0$};
        
\end{tikzpicture}
\caption{Scheme of (a) the boundary of the subset $U$ of the perturbation domain $\delta \V$ and (b) its details on a planar section. Blue (subindex $0$) and green (subindex $p$) represent the perturbed and unperturbed domains, respectively, whereas orange represents the generated boundary.}
\label{PertDomDetail}
\end{figure}
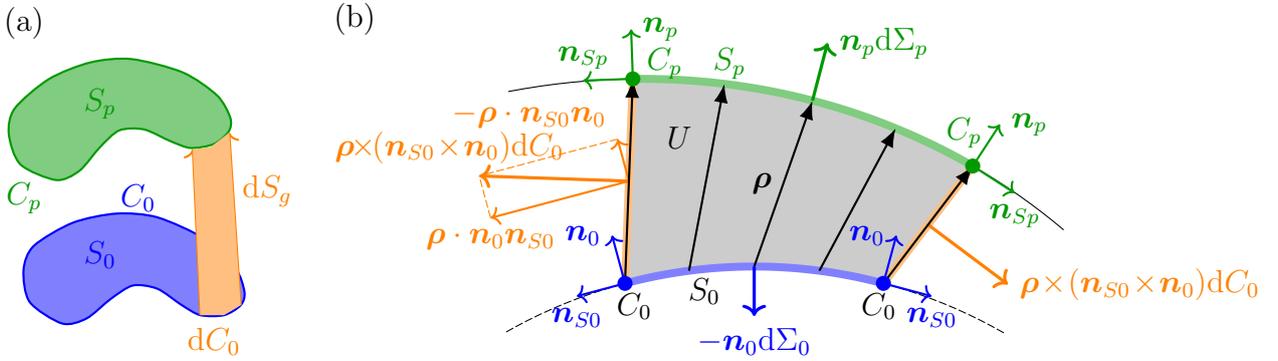

The perturbation domain $\delta \V$ is bounded by both the unperturbed and perturbed boundaries. It is convenient to consider the subset $U$ of the perturbation domain $\delta \V$ generated by sweeping a subset $S_0$ of the unperturbed boundary $\Sigma_0$ along an infinitesimal displacement $\vdn$ into a subset $S_p$ of the perturbed boundary $\Sigma_p$, as illustrated in \figref{PertDomDetail}. This subset domain $U$ is bounded by $S_0$, with differential vector $ - \vn_0 \dd \Sigma_0$, by $S_p$, with differential vector $ \vn_p \dd \Sigma_p$, and by the surface $S_g$ generated by sweeping the contour of $S_0$, denoted $C_0$, along the infinitesimal displacement $\vdn$, with differential vector $ \vdn \x \left( \vn_{S0} \x \vn_0 \right)  \dd C_0 $. %The components in $\vn_{S0}$ and $\vn_0$ of the differential vector of the generatrix can be obtained by developing the double vector product.

%we depict the details of a differential perturbation domain, which is bounded by the differential boundary on the unperturbed domain, with vector $ - \vn_0 \dd \Sigma_0$, lying on the unperturbed boundary, its perturbed counterpart with vector, $ \vn_p \dd \Sigma_p$, and the differential generatrix, with vector $  \vdn \x \left( \vn_{S0} \x \vn_0 \right)  \dd \Gamma_0 $, which is generated by sweeping a differential contour, see \figref{DomDecomp}(b), of the differential unperturbed surface $\dd \Gamma$ an amount $\vdn$, where $\vn_{S0}$ is the unperturbed counterpart of $\vnS$.  

Next, we explain the influence of the perturbation introduced by the infinitesimal displacement $\vdn$ on the integral over the perturbed domain and boundary, as well as on mixed boundary conditions and the exterior differential operator. In order to reduce the degrees of freedom of $\vdn$ to one, either the BALE method can be used, i.e.
\begin{align}\label{hbale}
\hdivS \gradS \vdn = h \, \vn \qquad \mbox{at } \Sigma_{0} \,, 
\end{align}
where $h$ represents the degree of freedom, or the displacement can be chosen to be in the normal direction, i.e. 
\begin{align}\label{rhon}
\vdn = \dn \, \vn \,.
\end{align} 
For this reason, we provide the DBP method for $\vdn$ fufilling either \eqref{hbale} or \eqref{rhon}. For the latter case we will see that many terms vanish.

\subsection{Perturbation of integrals}\label{PInt}

In this subsection, we express the integral over the perturbed domain and the perturbed boundary in terms of integrals over the unperturbed domain, boundary and contour. 

First, the integral of a generic quantity $\varphi$ over the perturbed domain $\V_p$ can be decomposed as the sum of the integral over the unperturbed domain plus perturbation terms. To obtain such a decomposition, the integral over a domain $\V$, such as the parameterisation is fixed $ \xi_i \in [\xi_{i}^{-}, \xi_{i}^{+}]$, of a quantity $\varphi$ is perturbed as 
\begin{align}\label{RTT0}
\intV{\varphi}{\V_p} = \intV{\varphi}{\V_0} +  \int_{\xi_1^{-}}^{\xi_1^{+}} \int_{\xi_2^{-}}^{\xi_2^{+}} \int_{\xi_3^{-}}^{\xi_3^{+}} \delta ( \varphi J)  \, \dd \xi_3 \, \dd \xi_2 \, \dd \xi_1 \,.
\end{align}
where the variation of the quantity is $\delta ( \varphi ) = \vdn \cdot \grad \varphi $ and of the Jacobian is $\delta ( J ) =J \,  \grad \cdot \vdn $, being $\vdn$ the infinitesimal displacement. Substituting these expressions in \eqref{RTT0} and using \eqref{VolGenStokes} leads to the Reynolds transport theorem
\begin{align}\label{RTT}
 \intV{\varphi}{\V_p} = \intV{ \left[ \varphi + \grad \cdot ( \vdn \varphi) \right]}{\V_0}   = \, \nonumber \\ 
                                = \intV{\varphi}{\V_0} + \intS{ \vn  \cdot \vdn \, \varphi}{\Sigma_0}
\end{align}
The unperturbed and perturbation contributions are respectively represented by the domain and boundary integral on the RHS. 

Second, the integral of a generic quantity $\varphi$ over the perturbed subset boundary $S_p$ can be decomposed as the sum of the integral over the unperturbed subset boundary plus perturbation terms. To obtain such a decomposition, the integral over a subset boundary $S$, such as the parameterisation is fixed $ \xi_i \in [\xi_{i}^{-}, \xi_{i}^{+}]$, of a quantity $\varphi$ is perturbed as
\begin{align}\label{BRTT0}
\intS{\varphi}{S_p} = \intS{\varphi}{S_0} +  \int_{\xi_1^{-}}^{\xi_1^{+}} \int_{\xi_2^{-}}^{\xi_2^{+}} \delta ( \varphi J_S)  \, \dd \xi_2 \, \dd \xi_1 \,.
\end{align}
where the variation of the quantity is $\delta (\varphi) = \vdn \cdot \grad \varphi$ and of the boundary Jacobian is $\delta ( J_S ) = J_S \, \gradS \cdot \vdn $, due to \cite[]{stone1990simple}. Substituting these expressions in \eqref{BRTT0} and using \eqref{equiv} and \eqref{VolGenStokes} leads to the boundary Reynolds transport theorem
\begin{align}
\label{intSdefS}
\int_{S_p} \varphi \, \dd \Sigma 
&= \int_{S_0} [ \varphi + \gradS \cdot (\vdn \varphi) + \vdn \cdot \vn\vn \cdot \grad \varphi ]  \, \dd \Sigma = \nonumber  \\
&= \int_{S_0} \varphi \, \dd \Sigma + \int_{C_{0}} \vnS \cdot \vdn \, \varphi \, \dd \Gamma + \int_{S_0}  \left[   \varphi \left( \vdn \cdot \vn \right) \divS \vn  +  \left( \vdn \cdot \vn \right) \left( \vn \cdot \grad \varphi \right) \right] \, \dd \Sigma  \,. 
\end{align}
The terms on the LHS and the first two terms on the RHS are analogous to \eqref{RTT}. However, it also appears an additional contribution due the non-euclidean variations of the boundary size and variations of the quantity in the direction out of the boundary.

Using \eqref{rhon} instead, eqs. \eqref{RTT} and \eqref{intSdefS} reduce to 
\begin{subequations}
\begin{align}
 \intV{\varphi}{\V_p} &= \intV{\varphi}{\V_0} + \intS{ \dn \, \varphi}{\Sigma_0} \,, \\
 \int_{S_p} \varphi \, \dd \Sigma &= \int_{S_0} \varphi \, \dd \Sigma 
+ \int_{S_0} \left[\dn \varphi  \divS \vn +  \dn \vn \cdot \grad \varphi \right]  \,  \dd \Sigma   \,. \label{Pert_int}
\end{align}
\end{subequations}
The equations \eqref{intSdefS} and \eqref{Pert_int} are also valid for $S_0 = \Sigma_0$ and $S_p = \Sigma_p$.

\subsection{Perturbation of mixed boundary conditions}\label{PMix}
   
In this subsection, we perturb the domain of a conservation PDE for a quantity $\phi$, subjected to a general boundary condition referred as a mixed or Robin boundary condition. For the sake of clarity, we consider the Poisson equation but the following procedure can be applied to any other conservation PDEs, as is done in \secref{Example1} to the Stokes equation. 
              
The aforementioned Poisson equation writes        
\begin{align}
\label{Poisson}
\div \grad \phi = \sigma \qquad \mbox{at } \V_p \,,
\end{align}     
where $\sigma=\sigma(\vx,\phi)$, together with the mixed boundary conditions
\begin{align}
\label{mixed}
\vn \cdot \grad \phi = c \phi + \gamma \qquad \mbox{at } \Sigma_p \,,
\end{align}     
where $\gamma = \gamma(\vx) $, and $c$ is a constant.

Considering that the domain $\V_p$ can be decomposed as in \figref{DomDecomp}, \eqref{Poisson} writes as
\begin{subequations}
\label{Poisson0}
\begin{align}
\div \grad \phi &= \sigma \qquad \mbox{at } \V_0 \,, \label{Poisson0a} \\
\div \grad \phi &= \sigma \qquad \mbox{at } \delta \V \,, \label{Poisson0b}
\end{align}
\end{subequations}     
For convenience and without loss of generality, \eqref{Poisson0b} is integrated over a subset $U$, shown in \figref{PertDomDetail}, of the perturbation domain $\delta \V$. 
In effect, using the Stokes theorem \eqref{VolGenStokes}, it writes as
\begin{align}
   \int_{S_p \cup S_0 \cup S_g} \vn \cdot \grad \phi \, \dd \Sigma 
=  \int_{U} \sigma \, \dd \V 
  \,, \qquad \forall \, S_0 \in \Sigma_0 \,.
\end{align}
where $\vn$ is the outer normal to the boundary $S_p$, $S_0$ and $S_g$, which up to first order in $\vdn$ rewrites as (see \figref{PertDomDetail})
\begin{align}
\label{Poisson1ini}
   \int_{S_p} \vn \cdot \grad \phi \, \dd \Sigma 
-  \int_{S_0} \vn \cdot \grad \phi \, \dd \Sigma 
+ \int_{C_0} [(\vn \x \vnS)  \x \vdn] \cdot \grad \phi \, \dd \Gamma
=
 \int_{S_0} \sigma \vdn \cdot \vn \, \dd \Sigma
  \,, \qquad \forall \, S_0 \in \Sigma_0 \,.
\end{align}
Using the boundary Stokes theorem \eqref{SurfGenStokes} and rearranging terms due to the flux through the perturbed boundary on the LHS, the flux through $S_p$ can be expressed in terms of the quantities defined at $S_0$ as
\begin{align}
\label{Poisson1}
 \int_{S_p} \vn \cdot \grad \phi \, \dd \Sigma 
= \int_{S_0} \{\vn \cdot \grad \phi + \hgradS \cdot [ \vdn \vn \cdot \grad \phi  - \vdn \cdot \vn \grad \phi ] + \sigma \vdn \cdot \vn \} \, \dd \Sigma \,, \qquad \forall \, S_0 \in \Sigma_0 \,. 
\end{align}          
The first term of the integral on the RHS represents the flux through $S_0$, the second and third ones represent the fluxes through the generatrix in the $\vn_{0}$ 
and $\vn_{S0}$ directions, respectively,  and the fourth one is the integral of the non-homogeneous term. Eq. \eqref{Poisson1} is the integral counterpart of \eqref{Poisson0b}.
%%Observe that the integrand of the LHS of \eqref{Poisson1} is the LHS of \eqref{mixed}. 

Integration of \eqref{mixed} over $S_p$, rewritten at $S_0$ with the use of \eqref{intSdefS} and \eqref{Poisson1}, leads for all $S_0 \in \Sigma_0$ to
\begin{align}
\label{mixedp}
\vn \cdot \grad \phi  &+ \hdivS \left( \vdn \vn \cdot \nabla \phi - \vdn \cdot \vn \grad \phi \right) + \sigma \vdn \cdot \vn   = \nonumber \\ 
&(c \phi + \gamma) + \divS \left[ \vdn \left( c \phi + \gamma \right) \right] + \left(\vdn \cdot \vn \right) \vn \cdot \grad (c \phi + \gamma)  \qquad \mbox{at } \Sigma_0 \,,
\end{align}     
which, using \eqref{rhon} instead, reduces to 
\begin{align}
\label{Pert_mix}
\vn \cdot \grad \phi   - \hgradS  \cdot \left( \dn  \grad \phi \right) + \dn  \sigma = 
\left[1 + \dn \left(\divS \vn \right) + \dn \vn \cdot \grad \right]\left(c \phi + \gamma\right) \qquad \mbox{at } \Sigma_0 \,.
\end{align}

 \subsection{Perturbation of the exterior differential operator}\label{AppD}

In this subsection, we perturb the exterior differential operator $\hgradS$. For this purpose, we first perturb the RHS of \eqref{SurfGenStokes},
\begin{alignat}{4}
\label{surfopdeflin00}
\int_{S_p} \hgradS \varphi \,\dd \Sigma =
\int_{C_p} \vnS \varphi \,\dd \Gamma =
\int_{C_0}  \vnS \varphi \,\dd \Gamma + \int_{C_0} [ \vnS \,\dd \Gamma \delta ( \varphi)   +   \delta( \vnS \dd \Gamma) \varphi  ] \,,
\end{alignat}
where $C_p$ is the contour $C_p \equiv \partial S_p$, and the variation of $\varphi$ can be written as
\begin{align}
\label{varvarphi}
\delta \! \left( \varphi \right)= \vdn \cdot \grad \varphi \,  
\end{align}
or $\delta \! \left( \varphi \right)= \vdn \cdot \gradS \varphi$ if $\varphi$ were only defined at the boundary. We will reproduce here the derivation for the variation of $\vnS \dd \Gamma$ given in \cite{Rivero2018aCorr} and which is analogous to the variation of a surface element \cite{batchelor2000introduction}.

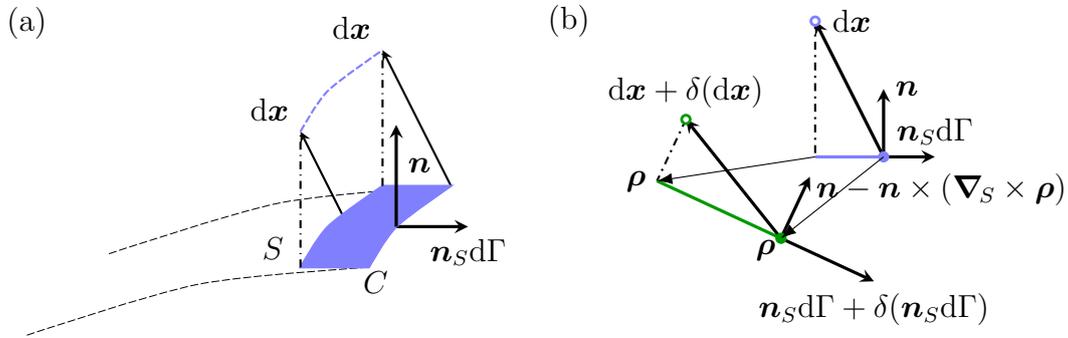
\begin{figure}[h!]
\centering
\begin{tikzpicture}[scale=.9,>=stealth]

%%% Sigma base
%\draw[very thin] plot [smooth] coordinates {(-5,-1) (-2.5,-.25) (0,0)};
\begin{scope}[shift={(.3,.3)}]
%\draw[very thin] plot [smooth] coordinates {(-5,-1) (-2.5,-.25) (0,0)};
\end{scope}
\begin{scope}[shift={(.6,.6)}]
\draw[dashed,very thin] plot [smooth] coordinates {(-5,-1) (-2.5,-.25) (0,0)};
\end{scope}
\begin{scope}[shift={(-.3,-.3)}]
%\draw[very thin] plot [smooth] coordinates {(-5,-1) (-2.5,-.25) (0,0)};
\end{scope}
\begin{scope}[shift={(-.6,-.6)}]
\draw[dashed,very thin] plot [smooth] coordinates {(-5,-1) (-2.5,-.25) (0,0)};
\end{scope}
%\draw[dashed,very thin] plot [smooth] coordinates {(-5.6,-1.6) (-4.4,-.4)};
\node[below] at (-2,0) {$S$};
\node[below] at (-.5,-.5) {$C$};

\draw [->,very thick] (-.2, 0)--(0.85, 0) node [below ] {$ \vGamma$} ;
%\draw [color1,ultra thick] (-.6,-.6)--(.6, .6);
\draw[color1, thick] plot [smooth] coordinates {(-.6,-.6)(-0.2, 0)(.6, .6)};

\draw [->,black,thick] (-.6,-.6)--(-1.6, 1.4) node[above left]{$\dd \vx$} ;
\draw [->,black,thick] (.6,.6)--(-.4, 2.6)  node[above left]{$\dd \vx$} ;
%\draw [color1,dashed, thick] (-1.6, 1.4)--(-.4, 2.6) ;
\draw[color1,dashed, thick] plot [smooth] coordinates {(-1.6, 1.4)(-1.2, 2)(-.4, 2.6)};

%\draw [color1,fill=color1,thick,opacity=0.2] (-.6,-.6)--(-1.6, -.6)--(-.4, .6)--(.6,.6) ;
\draw[color1,fill=color1,thick,opacity=0.2]  plot [smooth] coordinates {(-.6,-.6)(-0.2, 0)(.6, .6)}--
 plot [smooth] coordinates {(-.4,.6)(-1.2, 0)(-1.6, -.6)}--cycle;

\draw [dash dot,thick] (-1.6, -.6)--(-1.6, 1.4) ;
\draw [dash dot,thick] (-.4, .6)--(-.4, 2.6) ;

\draw [->,very thick] (-.2,0)--(-.2, 1.5) node[pos=.6,right] {$\vn$};
%\node[color1,below] at (-.5,-.7) {$\dd \Gamma$};

\node at (-5.6,3) {(a)};

\end{tikzpicture}
\begin{tikzpicture}[scale=.9,>=stealth]

%%% Sigma base
\draw[very thick,->] (0,0)--(-1,2) node[right] {$\, \dd \vx$} ;
\draw[very thick,->] (0,0)--(.75,0) node[above] {$ \vGamma$} ;
\draw[very thick,->] (0,0)--(0,1) node[right] {$\vn$} ;

\draw[thick,dash dot] (-1,2)--(-1,0);
\node at (0,0) [circle,fill=color1,inner sep=1.5pt]{};
\node at (-1,2) [circle,fill=color1,inner sep=1.5pt]{};
\node at (-1,2) [circle,fill=white,inner sep=.75pt]{};

\draw[color1,very thick,opacity=0.4] (0,0)--(-1,0);

%%%%%%%%%%%%%%%%%%%%%%%%%%%%

\begin{scope}[shift={(-1.5,-1.2)},rotate=-25]

%%% Sigma base
\draw[very thick,->] (0,0)--(-2,1) node[above] {$\dd \vx  + \delta \! \left( \dd \vx  \right) $} ;
\draw[very thick,->] (0,0)--(1.5,0) node[below] {$ \vGamma  + \delta \! \left( \vGamma  \right) $} ;
\draw[very thick,->] (0,0)--(0,1) node[right,pos=0.8] {$\vn-\vn \x (\rotS \vdn)$} ;

\draw[thick,dash dot] (-2,1)--(-2,0);
\node at (0,0) [circle,fill=green,inner sep=1.5pt]{};
\node at (-2,1) [circle,fill=green,inner sep=1.5pt]{};
\node at (-2,1) [circle,fill=white,inner sep=.75pt]{};

\draw[green,very thick,opacity=0.4] (0,0)--(-2,0);

\end{scope}

\draw[-{Latex[length=2mm]}] (0,0) -- (-1.5,-1.2) node[pos=1.15] {$\vdn$};
\draw[-{Latex[length=2mm]}] (-1,0) -- (-3.3126,-0.35476) node[pos=1,left] {$\vdn$};

\node at (-4.6,2) {(b)};

\end{tikzpicture}
\caption{Scheme of (a) the projection  on $S$ of the virtual boundary generated by sweeping the differential contour $\dd \Gamma$ along the virtual displacement $\dd \vx$, denoted $\dd \vx \cdot \vnS \dd \Gamma$, and (b) its variation represented before (blue) and after (red) the effect of the infinitesimal displacement $\vdn$.}
\label{SketchPert_MixDs}
\end{figure}

For this purpose, we consider the projection on the boundary $\Sigma$ of the virtual boundary, generated by sweeping the differential contour $\dd \Gamma$ along the virtual displacement $\dd \vx $, which is denoted by $\dd \vx \cdot \vnS \dd \Gamma$ and illustrated in \figref{SketchPert_MixDs}a. Its variation, illustrated in \figref{SketchPert_MixDs}b, can be written as
 \begin{align}
\label{vnS_dG0}
\delta ( \dd \vx \cdot \vnS \, \dd \Gamma) = \delta (\dd \vx) \cdot \vnS \dd   \Gamma + \dd \vx \cdot \delta ( \vnS \dd \Gamma) \,,
\end{align}
where, according to \cite{stone1990simple},
\begin{align}
\label{eq30}
\delta ( \dd \vx \cdot \vnS \dd \Gamma) =  \dd \vx \cdot \vnS \dd \Gamma \,  \divS \vdn \,.
\end{align}
The perturbation of the virtual vector $\dd \vx$ can be written as 
\begin{align}
\label{eq31}
\delta ( \dd \vx ) = \dd \vx \cdot \gradS \vdn + \dd \vx \cdot \vn \delta (\vn)  \,,
\end{align}
where the first term of the RHS represents the perturbation of the tangent component $\dd \vx \cdot \idS$ and the second term represents the perturbation of the normal component $ \left( \dd \vx \cdot \vn \right) \vn$ due to the rotation of the boundary, see \cite{weatherburn2016differential}
\begin{align}
\label{deltan}
\delta (\vn ) = - \vn \times \left( \rotS \vdn \right) = -\left( \gradS \vdn \right) \cdot \vn \,.
\end{align}
where the double vector product has been developed.
Introducing \eqref{eq30}, \eqref{eq31} and \eqref{deltan} in \eqref{vnS_dG0} which is valid for any $\dd \vx$, leads to
\begin{align}
\label{eq32}
\delta ( \vnS \dd \Gamma) &= 
\vnS \dd \Gamma \cdot [ \id \divS \vdn - (\gradS \vdn)^T + (\gradS \vdn) \cdot \vn\vn  ]\,.
\end{align}
An analytical derivation of \eqref{eq32} is provided in \ref{AppA}.

Thus, introducing \eqref{varvarphi} and \eqref{eq32} in \eqref{surfopdeflin00} and using \eqref{SurfGenStokes}, leads to
\begin{alignat}{4}
\label{surfopdefpert}
\intS{ \hgradS \varphi }{S_p} 
   = \intS{\hgradS \varphi}{S_0}+
   \intS{  \hdivS  \big{\{}  \big{[}  \id ( \divS \vdn )  -  (\gradS \vdn)^T + (\gradS \vdn) \cdot \vn\vn + \vdn \cdot \grad   \big{]} \varphi \big{\}}  }{S_0} \,,
\end{alignat}
which, using \eqref{rhon} instead, reduces to
\begin{alignat}{4}
\label{Pert_ext}
\intS{ \hgradS \varphi }{S_p}
   = \intS{  \hgradS  \varphi }{S_0} 
+ \intS{  \hdivS \big{\{}   \big{[}  \id (\dn \divS \vn  ) -  \dn \gradS \vn + (\gradS \dn) \vn  + \dn \vn \cdot \grad   \big{]} \varphi \big{\}} }{S_0} \,,
\end{alignat}
where $(\gradS \vn)^T = (\gradS \vn)$ and $ \vnS \cdot \vn =0$ have been used.

%After some algebraic manipulation, using the vector triple product and $ \vn \cdot \vnS = 0$, one obtains that the variation of the curve element is due to the dilatation of the contour and due to the rotation of the surface, respectively represented by the terms of the RHS of
%\begin{align}
%\vnS \dd \Gamma \cdot [ \id_S \dn \divS \vn -  (\gradS \dn \vn)^T + (\gradS \dn)\vn ]  = \vnS \dd \Gamma  \cdot ( \id \dn \divS \vn - \dn  \gradS \vn ) - \vnS \dd \Gamma \x [ \rotS (\dn \vn) ] \,.
%\end{align}

%%%%%%%%%%%%%%%%%%%%%
%%%%%%%%%%%%%%%%%%%%%
%%%%%%%%%%%%%%%%%%%%%

\section{Linearisation of partial differential equations at deformable domain using DBP} \label{Linearisation}

The linearisation of partial differential equations on deformable domains consists of two steps. The first step is the application of the DBP method in order to write the PDE at the unperturbed domain and its boundary as already done in \secref{DBP}. The second step is the asymptotic expansion of the variables into the PDE at the unperturbed domain and boundary as performed in this section. To this end, the infinitesimal displacement $\vdn$ and any variable $\varphi$ are expanded up to first order in $\epsilon\ll 1$ as 
\begin{subequations}
\label{expans}
\begin{align}
\vdn &\approx \epsilon \dn_1 \vn \,, \\
\varphi &\approx \varphi_0 + \epsilon \varphi_1 \,. \label{expansb}
\end{align}
\end{subequations} 
Substituting \eqref{expans} in the domain and boundary integrals \eqref{RTT} and \eqref{intSdefS} after application of the DBP method, leads up to first order in $\epsilon$ to
\begin{subequations}
%%%%\label{Pert_int}
\begin{align}
\int_{\V_p} \varphi \,  \dd \V &\approx \int_{\V_0} \varphi_0 \,  \dd \V + \epsilon \left ( \int_{\V_0} \varphi_1 \,  \dd \V + \int_{\Sigma_0} \varphi_0 \, \dn_1 \,  \dd \Sigma \right)  \,, \\
\int_{S_p} \varphi \, \dd \Sigma &\approx \int_{S_0} \varphi_0 \,  \dd \Sigma+ \epsilon  \int_{S_0} \left[ \varphi_1 + \left( \dn_1 \divS \vn  + \dn_1 \vn \cdot \grad \right) \varphi_0 \right]  \,  \dd \Sigma  \,.
\end{align}
\end{subequations}
 
For the mixed boundary conditions, the linearisation of the system \eqref{Poisson}, \eqref{mixed} is obtained by substituting \eqref{expans} into the system after the application of the DBP method, i.e. \eqref{Poisson0a} and \eqref{mixedp}. This leads to
\begin{subequations}
\label{Poisson0all}
\begin{align}
\grad \cdot \grad \phi_0 &= \sigma_0 \qquad \mbox{at } \mathcal{V}_0 \,,\\
 \vn     \cdot \grad \phi_0 &= c \varphi_0 + \gamma   \qquad \mbox{at } \Sigma_0 \,,
\end{align}  
\end{subequations}  
for the zeroth order in $\epsilon$, and
\begin{subequations}
\label{Poisson1all}
\begin{align}
\div \grad \phi_1 &= \sigma_{0,\phi} \phi_1 \qquad \mbox{at } \mathcal{V}_0 \,, \\
  \vn \cdot \grad \phi_1  - \hgradS \cdot \left( \dn_1 \grad \phi_0 \right) + \dn_1 \sigma_0 &=  c \phi_1 +  ( \dn_1 \divS \vn + \dn_1 \vn \cdot \grad ) (c \phi_0 + \gamma)
   \qquad \mbox{at } \Sigma_0 \,,
\end{align}        
\end{subequations}
for the first order in $\epsilon$, where we have Tayloer expanded $\sigma \approx \sigma_0+ \epsilon \sigma_{0,\phi} \phi_1  $, with $\sigma_0 = \sigma(\phi_0,\vx)$ and $\sigma_{0,\phi} = \partial_{\phi} \sigma(\phi_0,\vx)$.

For the linearisation of the exterior differential operator $\hgradS \phi$, \eqref{expans} is substituted in its expression after the DBP method \eqref{surfopdefpert} which leads, up to first order, to
\begin{alignat}{4}
\label{Pert_ext01}
\intS{ \hgradS &\varphi }{S_p}
   = \intS{  \hgradS  \varphi_0 }{S_0} \nonumber\\
&+ \epsilon \intS{ \big{\{} \hdivS \big{[}  (\dn_1 \divS \vn  ) \id -  \dn_1 \gradS \vn + (\gradS \dn_1) \vn  + \dn_1 \vn \cdot \grad   \big{]}  \varphi_0 + \hgradS \varphi_1 \big{\}} }{S_0} \,,
\end{alignat}

\section{Capillary migration of bubbles in microchannels} \label{Example1}

In this section, we exemplify the previous methods by their application to a problem with deformable domain whose final domain is part of the solution of the problem. We have used these methods to solve a problem of interest such as the transverse migration force experienced by a deformable bubble flowing inside a microchannel \cite[]{Rivero2018a}. In our previous work, we have carried out singular asymptotic expansion on the $\Ca$ number around zero for which the bubble shape is spherical. In the present work, for the sake of clarity and exposition of the method, we apply this method to a related problem which requires a regular asymptotic expansion, instead of singular, i.e. for a finite $\Ca$ number for which the shape of the bubble is not known a priori. 

\begin{figure}[h!]
\centering
\begin{tikzpicture}[x=.08\textwidth,y=.08\textwidth]

%\draw[dashed]  
%plot [domain = 0:360, samples = 40, variable = \i]  ({1.5+(1+.0*\i/100)*cos \i}, {2.5+(1+.0*\i/100)*sin \i}, 0); 

%\draw[fill=cyan!25, dashdot] (-1.5,2.5) circle (.8);     \node at (-1.5,2.5)   {$\VR$};
\begin{scope}[shift={(-1.5,2.5)}]
\draw[fill=cyan!25, dashdot,rounded corners=10pt] (-.7,.7)--(-.7,-.7)--(.7,-.7)--(.7,.7)--cycle;     \node at (0,0)   {$\vX \in \V_{Mat}$};   \node at (.85,.85)   {$\Sigma_{Mat}$};
\node at (3,-1)   {$0$-th order};
\node at (6,-1)   {$0$-th + $1$-st order};
\end{scope}

\draw[fill=cyan!25, dashed] (1.5,2.5) circle (.8);    \node at (1.5,2.5)   {\small $\vx \!\in\! \V_{Spa} \!\equiv\! \V_{0}$}; \node at (2.5,3.3)   {$\Sigma_0 \!\equiv\! \Sigma_{Spa}$};
\draw[->,thick] (-.5,2.5)--(.5,2.5) node[midway,above] {$\vq$}  node[midway,below] {BALE};

%\draw[->,thick] (-.7,1.7)--(1.3,.7) node[pos=.3,anchor=south west] {$\vq+ \vdn$};
%\draw[->,thick] (1.8,1.6)--(2.0,.85) node[midway,right] {$\vdn$} node[midway,left] {DBP} ;

\begin{scope}[shift={(4.5,2.5)}]   
\draw[fill=cyan!25,rotate=45] (0,0) ellipse  (.9 and .7)  ;        
\node at (0,0) {$\V_p$} ;
\node at (.85,.8) {$\Sigma_p$} ;
\end{scope}   

\begin{scope}[shift={(3,0)}]   
\draw[->,thick] (-.5,2.5)--(.5,2.5) node[midway,above] {$\vdn$}  node[midway,below] {DBP};
\end{scope}

%\draw[fill=cyan!25,cyan!25,rotate=45, even odd rule] (0,0) circle (.8)  (0,0) ellipse  (.9 and .7)  ;                 
%\draw[rotate=45] (0,0) ellipse  (.9 and .7);
%\draw[dashed] (0,0) circle  (.8);
%\node at (.5,.5)   {$+$};
%\node at (-.5,-.5)   {$+$};
%\node at (.4,-.4)   {$-$};
%\node at (-.4,.4)   {$-$};
%     
%\node at (-2,0)   {$\V_0$};
%\node at (-1,0)   {$\cup$};
%\node at (0,0)   {$ \delta \V $};
%\node at (1,0)   {$\equiv$};
%\node at (2,0)   {$\V_p $};            
%
%\node at (-1.2,.8)   {$\Sigma_0$};
%\node at (2.8,.8)   {$\Sigma_p$};
%
%\node at (0,-1.5)   {};

\end{tikzpicture}
\caption{Scheme of domain deformation (BALE) and its perturbation (DBP) used for the problem of capillary migration.}
\label{Sketch1}
\end{figure}
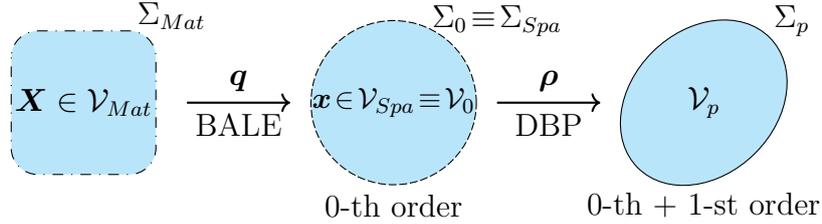

In this example, we have linearised the system of equations using the linearisation method outlined in \secref{Linearisation}, which is based on the DBP method presented in \secref{DBP}. This results in two system of equations: (i) a nonlinear one, which is the same as the original system of equation, defined in an unknown domain, which can be modelled with the use of the BALE method presented in \secref{BALE} and (ii) a linear system which is defined in the previous domain containing the terms rising from the use of the DBP method. In this case, the spatial domain $\V_{Spa}$ used in the BALE method coincides with the undeformed domain $\V_0$ used in the DBP method, as illustrated in \figref{Sketch1}. It is worth mentioning that both methods are independent and can be used combined as in this case, or separately as in our previous work \cite[]{Rivero2018a}, where, in the one hand, the nonlinear equation is solved, and in the other hand, the linearisation is done for a case in which the unperturbed domain is a priori known.

For the sake of completeness, and prior to that, we briefly model the aforementioned system in strong formulation already available in \cite{Rivero2018a}. Then, we provide with the linearisation in strong formulation, and we finally write the weak formulation of both.
 
\subsection{Governing equations}\label{zeroth}

\begin{figure}[h!]
\centering
\begin{tikzpicture}[x=.5\textwidth,y=.5\textwidth]

\definecolor{mycolor1}{rgb}{0.1,0.1,.1}%
\definecolor{mycolor2}{rgb}{.4,.4,.4}%

	\draw[fill=mycolor1!40,draw=black,opacity=1] (-.5,-.2) rectangle (.5,.2);

	\node at (-.45,.1) {$  \Sigma_{out}$} ;
	\node at (.54,.1) {$  \Sigma_{in}$} ;
	\node at (0.1,+0.23) {$  \Sigma_{W}$} ;
	\node at (0.1,-0.17) {$  \Sigma_{W}$} ;

	\draw[fill=mycolor2!40,draw=black,opacity=1] plot [domain = 0:360, samples = 80, variable = \i] ({.08*(1+0*cos (45+2*\i))*cos \i}, {.07+.08*(1+0*cos (45+2*\i))*sin \i}, 0) ;

	\draw[thick, |->] (0,0) -- (0,0.07) node[midway,right] {$\veps$};
	
	\draw[thick, |->] (.34,.12) -- (.34,.04) node[midway,right] {$\vf$};
	
	\draw[black,dash dot] (-.52,0) -- (.52,0) ;
	
	\node at (0.11,0.11) {$  \Sigma_{B}$} ;
	
	\node at (-0.2,-0.12) {$  \mathcal{V}$} ;
	\node at (-0.02,0.11) {$  \mathcal{V}_{B}$} ;

	\draw[<->] (-.5,-.23) -- (.5,-.23) node[midway,below]{$L$} ;
	\draw[<->] (-.53,-.2) -- (-.53,.2) node[midway,left]{$1$} ;

	\draw[thick,->] (.3,-.12) -- (.3,-.04) node[left]{$\ve_y$} ;	
	\draw[thick,->] (.3,-.12) -- (.38,-.12) node[below]{$\ve_x$} ;	
	
%	\draw[fill=mycolor1!40,draw=mycolor1!40,opacity=.5] (-.5,-.2) rectangle (.5,.2)
%	plot [domain = 0:360, samples = 80, variable = \i] ({.15*(1+.1*cos (45+2*\i))*cos \i}, {.1+.15*(1+.1*cos (45+2*\i))*sin \i}, 0) ;

%	\draw[fill=mycolor2!40,draw=mycolor2!40,opacity=.5]  
%	(-.5,-.3) rectangle (.5,.3)
%	plot [domain = 0:360, samples = 80, variable = \i]  (.15*cos \i, .1+.15*sin \i, 0);% node[left] {$\Sigma_0$}  ;
%
%
%	\draw[line width=.2] plot [domain = 0:360, samples = 80, variable = \i] 
%                   ({.15*(1+.1*cos (45+2*\i))*cos \i}, {.1+.15*(1+.1*cos (45+2*\i))*sin \i}, 0) ;
%	
%	\draw[line width=.2,dashed]  plot [domain = 0:360, samples = 80, variable = \i] 
%                    (.15*cos \i, .1+.15*sin \i, 0);% node[left] {$\Sigma_0$}  ;
%
%	\draw[line width=.2,fill=mycolor1!40,draw=mycolor1!40,dashed,opacity=.5]  	(.51,-.02) rectangle (.65,.08);
%	\node at (.58,.03) {$  \mathcal{V}$} ; 
%
%	\draw[line width=.2] (.51,-.13) rectangle (.65,-.03);
%	\node at (.58,-.08) {$  \Sigma_B$} ; 
%
%	\draw[line width=.2,fill=mycolor2!40,draw=mycolor2!40,dashed,opacity=.5]  	(.51,.2) rectangle (.65,.3);
%	\node at (.58,.25) {$  \mathcal{V}_0$} ; 
%
%	\draw[line width=.2,dashed] (.51,.09) rectangle (.65,.19);
%	\node at (.58,.14) {$  \Sigma_{B_0}$} ; 
%
%\node[anchor=south east] at (-.5,.3) {(a)};

\end{tikzpicture}
\caption{Sketch of the geometry including domains $\V$ and $\V_B$ as well as the contours $\Sigma_B$, $\Sigma_{W}$, $\Sigma_{in}$ and $\Sigma_{out}$.}
\label{Geo}
\end{figure}
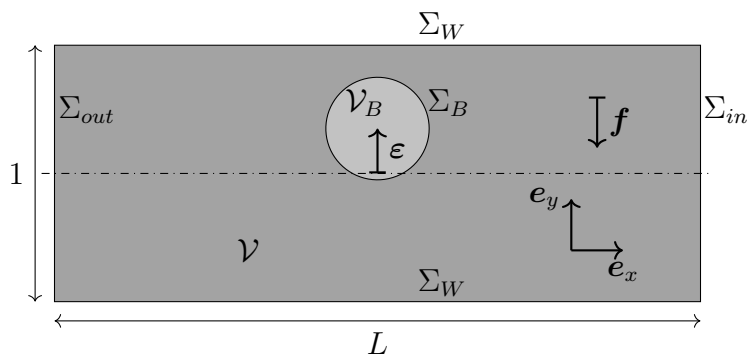

In this subsection, we introduce the dimensionless equations governing the steady dynamics of a periodic train of bubbles in microchannels which are defined on the domain shown in \figref{Geo}. In this physical situation, we consider a bubble of volume $\V_B$ and pressure of the gas $p_G$ located at $\veps$ within a periodic domain $\V$ of period $L$ in the longitudinal direction $x$. The bubble travels at a constant velocity $V$ in the $x$ direction and the reference frame moves attached to the bubble. The unitary flow field produces a Poiseuille pressure drop $L \partial_x p_P$ along a period $L$ modified by the presence of the bubble $\Delta p$ as well as an hydrodynamic force which is in equilibrium at $\veps$ with the buoyancy due to an uniform force exerted on the liquid $\vf$ in the transverse direction, i.e. $\vf \cdot \ve_x = 0$. In this work, and for the sake of simplicity, we consider 2D geometries which are periodic in the longitudinal direction $x$. Therefore domains are areas and boundaries are lines. The domain consists in a rectangle $\V$ of width $L$ and height $1$ with a hole $\V_B$ located at $\veps$ which is centred in the horizontal direction and with an off-centred position such as $\veps  = \varepsilon \ve_y $ from the centreline of the microchannel. The upper and lower sides of the rectangle represent the wall of the microchannel $\Sigma_W$ and the left and right sides, denoted $\Sigma_{out}$ and $\Sigma_{in}$, represent two cross sections of the microchannel, whereas the contour of the hole, $\Sigma_B$ represents the interface of the bubble. For more details on the modelisation, the reader is referred to our previous work \cite[]{Rivero2018a}.

The flow is governed by the Stokes equations written in dimensionless form,
\begin{subequations}
\label{Stokes}
\begin{align}
\div \vv &= 0  \qquad \mbox{at $\V$} \,,\\  
\div \stress &= \vzero \qquad \mbox{at $\V$} \,,
\end{align}
\end{subequations}
where $\vv$ and $\stress = - p \id + \grad \vv + (\grad \vv)^T$ are the velocity and the reduced stress tensor, and $p$ is the reduced pressure with a reference at a point $\vx_p$, 
\begin{align}\label{pref}
p = 0  \qquad \mbox{at $\vx_p$} \,.
\end{align}
The walls velocity is
\begin{align}\label{Wall}
\vv  = - V \ve_x  \qquad \mbox{at $\Sigma_W$} \,.
\end{align}
The impermeability condition writes
\begin{align}
\label{imp}
 \vn \cdot \vv &= 0   \qquad \mbox{at $\Sigma_B$} \,,
 \end{align}
and the stress balance is governed by the Young-Laplace equation
\begin{align} 
 \label{YL}
 \vn \cdot \llbracket \stress \rrbracket &= \frac{1}{Ca} \hgradS 1  \qquad \mbox{at $\Sigma_B$} \,,
\end{align}
where $\llbracket \stress \rrbracket = \stress + \id [p_G - \vf \cdot (\vx - \veps)]$ is the stress jump which contains the hydrostatic pressure $\vf \cdot (\vx - \veps)$ and $\hgradS 1$ is, according to \eqref{equiv}, the mean curvature vector of the boundary, i.e. $\hgradS 1 = \gradS \cdot \idS$. In addition, periodicity conditions in the longitudinal direction require
\begin{subequations}
\label{periodicity}
\begin{alignat}{4}
 p(\vx) &= p(\vx+L\ve_x) + \Delta p - L \partial_x p_P &\qquad& \mbox{at $\Sigma_{in}$} \label{pjump} \,,\\
\label{perv} \vv(\vx) &= \vv(\vx+L\ve_x) &\qquad& \mbox{at $\Sigma_{in}$} \,,\\
\label{perdnv} \vn \cdot \grad \vv(\vx) &= \vn \cdot \grad \vv(\vx+L\ve_x) &\qquad& \mbox{at $\Sigma_{in}$} \,,
\end{alignat}
\end{subequations}
where the 2D Poiseuille pressure drop is $\partial_x p_P = -12 $ and produces an unitary flow rate, whence
\begin{align}\label{Flow}
0 &= \int_{\Sigma_{in}} (\vv \cdot \ve_x + V - 1 )   \, \dd \Sigma    \,. 
\end{align}
We consider that the system is at equilibrium when the longitudinal force exerted on the bubble vanishes, i.e.
\begin{align}\label{Drag}
0= \vf \cdot \ve_x    \,.
\end{align}
The size of the domain occupied by the bubble, $\V_B$, and the geometric centre of the bubble, $\veps$, can be defined as 
\begin{subequations}\label{Surface}
\begin{align} 
\V_B &= \int_{\V_B}   \, \dd \V    \,, \\
\V_B \veps &= \int_{\V_B} \vx   \, \dd \V    \,, 
\end{align}
\end{subequations}
where the domain occupied by the bubble is a priori unknown. The system of equations \eqref{Stokes}-\eqref{Surface} is time dependent with domain variables $p$ and $\vv$, together with the surface variable $\dn$ and the global variables $\vf$, $V$, $\Delta p$ and $p_G$, whereas the geometry is unknown. The values of $\Ca$, $\veps$ and $\V_B$ are known. 

%%%%%%%%%%%
%%%%%%%%%%%
%%%%%%%%%%%
%%%%%%%%%%%
%%%%%%%%%%%
%%%%%%%%%%%

\subsection{Linearisation}\label{lin}

In order to perturb the system \eqref{Stokes}-\eqref{Surface} using the proposed method, the DBP method is first applied for $\vdn= \dn \vn$ and the variables are next expanded as \eqref{expans}. 

First, using the DBP method, the Stokes equations \eqref{Stokes} as well as the impermeability \eqref{imp} and stress balance \eqref{YL} boundary conditions write, using \eqref{Pert_mix} and \eqref{Pert_ext}, at the unperturbed domain as
\begin{subequations}
\label{StokesV0}
\begin{alignat}{5}
\div \vv       &=& 0         &\qquad& \mbox{at $\V_0$} \,,\\
\div \stress &=& \vzero &\qquad& \mbox{at $\V_0$} \,,  
\end{alignat}
\end{subequations}
and at the unperturbed boundary as 
\begin{subequations}
\begin{align}
\vn \cdot  \vv  -\hgradS \cdot (\dn \vv)  &= 0  \qquad \mbox{at $\Sigma_{B0}$}\,, \label{impV0} \\
\vn\cdot  \llbracket \stress \rrbracket -\hgradS \cdot ( \dn \llbracket \stress \rrbracket) - \dn \vf  
 &=  \frac{1}{Ca}  \hdivS [(1 + \dn \divS \vn) \id - \dn \gradS  \vn + (\gradS \dn)\vn ]  \qquad \mbox{at $\Sigma_{B0}$} \,,  \label{YLV0}
\end{align}
\end{subequations}
where the integral of $\grad \cdot \{ \id [p_G  - \vf \cdot (\vx-\veps)] \} + \vf = \vzero $ over $U$ has also been used. The boundaries $\Sigma_W$,  $\Sigma_{in}$ and $\Sigma_{out}$ are not deformed and thus \eqref{Wall}, \eqref{periodicity} and \eqref{Flow} are not affected, as well as the global equations \eqref{Drag}. The domain of \eqref{Surface} is also perturbed and then writes, using \eqref{RTT}, as 
\begin{subequations}
\label{Surface0}
\begin{align}
\V_B &= \intV{}{\V_{B0}} + \intS{ \dn}{\Sigma_{B0}}  \,, \\
\V_B \veps &=  \intV{\vx}{\V_{B0}} + \intS{ \vx \dn}{\Sigma_{B0}}  \,.
\end{align}
\end{subequations} 

%%% After decomposing the perturbed volume, the full system consists of \eqref{StokesV0}, \eqref{pref}, \eqref{Wall}, \eqref{imp0}, \eqref{YL0}, \eqref{periodicity}, \eqref{Flow}, \eqref{Drag} and \eqref{Surface0}. 

%In this work, we linearise $p \approx p_0 + \varepsilon p_1$, $\vv\approx \ve_x (u_0 + \varepsilon u_1) + \ve_y (v_0 + \varepsilon v_1) $, $\vf \approx \ve_y (f_0 + \varepsilon f_1)$, $V \approx V_0 + \varepsilon V_1$, $\delta p \approx \delta p_0 + \varepsilon \delta p_1$, $\dn \approx \varepsilon \dn_1)$ and $p_G \approx (p_{G0} + \varepsilon p_{G1})$. 
 
Second, expanding the variables in the previous system as \eqref{expans} where $\varphi$ is any variable $p$, $\vv$, $V$, $\Delta p$, $\vf$ and $p_G$ as well as $\veps = (\varepsilon + \epsilon) \ve_y $ and $\dn= \epsilon \dn_1 $ leads to the system described in what follows. The Stokes equations \eqref{StokesV0} leads to
\begin{subequations}
\label{Stokes1}
\begin{alignat}{8}
\div \vv_0       &=  0\,,        &\qquad&\div \vv_1        &=  0        &&\qquad& \mbox{at $\V_0$} \,, \label{Stokes1a} \\
\div \stress_0 &= \vzero\,, &\qquad& \div \stress_1 &= \vzero  &&\qquad& \mbox{at $\V_0$} \,, \label{Stokes1b}
\end{alignat}
\end{subequations}
together with the pressure reference at a point $\vx_p$ \eqref{pref}, 
%\begin{subequations}
\begin{align}
\label{pref1}
p_0 = 0 \,, \qquad p_1 = 0  \qquad \mbox{at $\vx_p$} \,,
\end{align}
%\end{subequations}
vanishing velocity of the wall \eqref{Wall}
%\begin{subequations}
\begin{align}
\label{Wall1}
\vv_0 + V_0 \ve_x = \vzero \,, \qquad  \vv_1 + V_1 \ve_x = \vzero  \qquad \mbox{at $\Sigma_W$} \,,
\end{align}
%\end{subequations}
impermeability condition \eqref{impV0},
%\begin{subequations}
\begin{alignat}{4}
\label{imp1}
 \vn \cdot \vv_0 = \vzero \,,\qquad \vn \cdot \vv_1 - \hdivS \left( \dn_1 \vv_0 \right)  &= \vzero    &\qquad& \mbox{at $\Sigma_{B0}$} \,,
 \end{alignat}
%\end{subequations}
stress balance \eqref{YLV0},
 \begin{align} \label{YL1}
\vn \cdot \llbracket \stress_0 \rrbracket &= \frac{1}{\Ca} \hgradS 1 \,, \nonumber \\
\vn \cdot \llbracket \stress_1 \rrbracket - \hgradS \cdot ( \dn_1 \llbracket \stress_0 \rrbracket) - \dn_1 \vf_0  &= \frac{1}{\Ca} \hdivS [ \id \dn_1 \divS \vn - \dn_1 (\gradS  \vn) + (\gradS \dn_1) \vn  ]   \qquad \mbox{at $\Sigma_{B0}$}
 \end{align}
where $\llbracket \stress_i \rrbracket = - p_i \id + \grad \vv_i + (\grad \vv_i)^T + \id [p_{G,i} - \vf_i \cdot (\vx - \veps)]$, the periodicity condition \eqref{periodicity},
\begin{subequations}
\begin{alignat}{8}
p_0(\vx) &= p_0(\vx \!+\!L\ve_x) - \Delta p_0 - L \partial_x p_P &&\quad&  p_1(\vx) &= p_1(\vx \!+\! L\ve_x) - \Delta p_1  &\quad& \mbox{at $\Sigma_{in}$} \label{pjump1} \,,\\
\label{perv1}  \vv_0(\vx) &= \vv_0(\vx \!+\!L\ve_x) &&\quad& \vv_1(\vx) &= \vv_1(\vx\!+\!L\ve_x)  &\quad& \mbox{at $\Sigma_{in}$} \,,  \\
\label{perdnv1} \ngrad \vv_0(\vx) &= \ngrad \vv_0(\vx\!+\!L\ve_x)  &&\quad&  \ngrad \vv_1(\vx) &= \ngrad \vv_1(\vx\!+\!L\ve_x)  \vzero &\quad& \mbox{at $\Sigma_{in}$} \,, 
\end{alignat}
\end{subequations}
average flow rate \eqref{Flow},
\begin{align}\label{Flow1}
 \intS{ (\vv_0 \cdot \ve_x + V_0 -1 )}{\Sigma_{in}} = 0 \,,\qquad
 \intS{ (\vv_1 \cdot \ve_x + V_1 )}{\Sigma_{in}} = 0  \,, 
\end{align}
equilibrium of the bubble in longitudinal direction \eqref{Drag},
\begin{align}
\label{Drag1}
 \vf_0 \cdot \ve_x = 0    \,,\qquad
 \vf_1 \cdot \ve_x = 0    \,,
\end{align}
and definition of the bubble size and position \eqref{Surface},
\begin{subequations}
\label{Surface1}
\begin{alignat}{6} 
\intS{ }{\V_{B0}}    &= \V_{B0}   \,,&&\qquad&
\intS{ \dn_1 }{\Sigma_{B0}}    &= 0   \,,  \\
\intS{\vx}{\V_{B0}}    &= \V_{B} \varepsilon \ve_y    \,,&&\qquad&
\intS{ \vx \dn_1}{\Sigma_{B0}}    &= \V_{B} \ve_y   \,.
\end{alignat}
\end{subequations}

\subsection{Weak formulation}\label{WF}
 
In this subsection, we write the weak form of the components of the previous systems of partial differential equations \eqref{Stokes1}-\eqref{Surface1}, being hereafter $x \equiv x_1$ and  $y \equiv x_2$ and variables with tildes representing the test functions of the variables without tildes. In the system \eqref{Stokes1}-\eqref{Surface1}, there are two equations per line, the first one corresponding to the zeroth order and the second one to the first order. Thus, in what follows \eqref{Stokes1}-\eqref{Surface1} refer either to the zeroth or to the first order accordingly.

Concerning the zeroth order equations, the vectors can be written in the laboratory frame as $\vv_0 = u_{0x} \ve_x + u_{0y} \ve_y $, $\vq = q_x \ve_x + q_y \ve_y $ and $\vf_{\!0} = f_0 \ve_y$, whence \eqref{Drag1} is automatically fulfilled. Multiplying the continuity equation \eqref{Stokes1a} by $\test{p}_0$, and imposing the pressure reference \eqref{pref1} using the Lagrange multiplier $\test{p}_{\rm{0ref}}$ which is a Dirac delta of unknown amplitude located at $\vx_{\rm{ref}}$, leads to 
\begin{align}
\label{pp0}
0 =& \intV{ \Big{[}\test{p}_0 ( \partial_x u_{0x} \!+\! \partial_y u_{0y} )    + \test{p}_{0 \rm{ref}} {p}_0 + {p}_{0\rm{ref}} \test{p}_0  \Big{]} }{\V_0}  
\end{align}
Multiplying the $x$ and $y$ components of the Stokes equations \eqref{Stokes1b} by $\test{u}_{0x}$ and $\test{u}_{0y}$ and integrating over the domain $\V_0$, as well as using the reciprocal theorems \eqref{prop3} and imposing the Dirichlet boundary conditions \eqref{Wall1} at $\Sigma_W$ using the Lagrange multipliers $\tau_{0nx}$ and $\tau_{0ny}$, \eqref{YL1} at $\Sigma_{B0}$ and \eqref{pjump1} at $\Sigma_{in}$, leads to
\begin{align}
\label{w0}
0 =& \intV{ \Big{[} - \partial_x \test{u}_{0x} (2 \partial_x u_{0x} \!-\! p_0 ) - \partial_y \test{u}_{0y} (2 \partial_y u_{0y}\!-\!p_0) - (\partial_y \test{u}_{0x} \!+\! \partial_x \test{u}_{0y})(\partial_y u_{0x} \!+\! \partial_x u_{0y}) \Big{]}  }{\V_0}  \nonumber\\ 
    +& \intS{ \Big{[} \test{u}_{0x} \tau_{0nx} + \test{u}_{0y} \tau_{0ny} + \test{\tau}_{0nx} (u_{0x}\!+\!V_0) + \test{\tau}_{0ny} u_{0y} \Big{]} }{\Sigma_W} \nonumber\\ 
    +& \intS{ \Big{[} (\test{u}_{0x} n_x + \test{u}_{0y} n_y) (-p_{0G} \!+\! f_{0} y) - \frac{1}{\Ca} ( \partial_{Sx} \test{u}_{0x} \!+\! \partial_{Sy} \test{u}_{0y}) \Big{]} }{\Sigma_{B0}} \nonumber\\ 
    +& \intS{ \test{u}_{0x} (\Delta p_0 +12  L )  }{\Sigma_{in}} \,,
 \end{align}
where variables $\tau_{0nx}$ and $\tau_{0ny}$ represent the stresses exerted on the wall in the $x$ and $y$ directions, respectively. The functions $u_{0x}$ and $u_{0y}$ are periodic, thus fulfilling the periodicity equations \eqref{perv1} and \eqref{perdnv1}. Multipliying the $x$ and $y$ components of \eqref{BALEeq} by $\test{q}_{x}$ and $\test{q}_{y}$, respectively, as well as the impermeability condition \eqref{imp1} by $\test{g}$ and integrating over the boundary $\Sigma_{B0}$ leads to 
 \begin{align}
0= \intS{  [ \partial_{Sx} \test{q}_{x} \partial_{Sx} {q}_{x} \!+\!  
                     \partial_{Sy} \test{q}_{x} \partial_{Sy} {q}_{x} +  
                     \partial_{Sx} \test{q}_{y} \partial_{Sx} {q}_{y} \!+\!
                     \partial_{Sy} \test{q}_{y} \partial_{Sy} {q}_{y} + \nonumber\\
                     g(\test{q}_x n_x \!+\! \test{q}_y n_y ) + \test{g} (u_x n_x \!+\! u_y n_y) ]}{\Sigma_{B0}} \,.
 \end{align} 
In addition, the global equations \eqref{Flow1} and \eqref{Surface1} as well as the partial differential equation \eqref{ALE} for the change of variable $\vx-\vX$ close the system of equations. 

Concerning the first order equations, the vectors can be written in the laboratory frame as $\vv_1 = u_{1x} \ve_x + u_{1y} \ve_y $, and $\vf_1 = f_1 \ve_y$, where \eqref{Drag1} is automatically fulfilled. Multiplying the continuity equation \eqref{Stokes1a} by $\test{p}_1$, and imposing the pressure reference \eqref{pref1} using the Lagrange multiplier $\test{p}_{\rm{0ref}}$ which is a Dirac delta of unknown amplitude located at $\vx_{\rm{ref}}$, leads to  
\begin{align}
\label{pp1}
0 = \intV{  \Big{[}\test{p}_{1} ( \partial_x u_{1x} \!+\! \partial_y u_{1y} )  + \test{p}_{1\rm{ref}} {p}_1 + {p}_{1\rm{ref}} \test{p}_1 \, \Big{]}  }{\V_0}  
\end{align}
Multiplying the $x$ and $y$ components of the Stokes equations \eqref{Stokes1b} by $\test{u}_{1x}$ and $\test{u}_{1y}$ and integrating over the domain $\V_0$, as well as using the reciprocal theorems \eqref{prop3} and imposing the Dirichlet boundary conditions \eqref{Wall1} at $\Sigma_W$ using the Lagrange multipliers $\tau_{1nx}$ and $\tau_{1ny}$, \eqref{YL1} at $\Sigma_{B0}$ and \eqref{pjump1} at $\Sigma_{in}$, leads to
\begin{align}
0 =& \intV{  \Big{[} - \partial_x \test{u}_{1x} (2 \partial_x  u_{1x} \!-\! p_1) - \partial_y \test{u}_{1y} (2 \partial_y v_{1y}\!-\!p_1) - (\partial_y  \test{u}_{1x} \!+\! \partial_x \test{u}_{1y})(\partial_y  {u}_{1x} \!+\! \partial_x {u}_{1y}) \Big{]}  }{\V_0}  \nonumber\\ 
    +& \intS{ \Big{[} \test{u}_{1x} \tau_{1nx} + \test{u}_{1y} \tau_{1ny} + \test{\tau}_{1nx} (u_{1x} \!+\! V_{1}) + \test{\tau}_{1ny} u_{1y} \Big{]} }{\Sigma_W}  \nonumber\\ 
    +& \intS{ \Big{\{} -\dn_1 [ \partial_{Sx} \test{u}_{1x} ( 2 \partial_{x} u_{0x} \!-\! p_0 \!+\! p_{0G} \!-\! f_{0} y) + \partial_{Sy} \test{u}_{1y} ( 2 \partial_{y} u_{0y}\!-\!p_0\!+\!p_{0G} \!-\! f_0 y)  \nonumber\\ 
    &\qquad + ( \partial_{Sy} \test{u}_{1x} \!+\! \partial_{Sx} \test{u}_{1y}) (\partial_{y} {u}_{0x} \!+\! \partial_{x} {u}_{0y})  ] + (\test{u}_{1x} n_x \!+\! \test{u}_{1y} n_y) (-p_{1G} \!+\! f_{1} y) + \dn_1  \test{u}_{1y}  f_0  
     \nonumber\\ 
    &\qquad- \frac{1}{\Ca} \, \dn_1  \left(\partial_{Sx} \test{u}_{1x} \!+\! \partial_{Sy} \test{u}_{1y})  (\partial_{Sx} n_x + \partial_{Sy} n_y \right)    \nonumber\\ 
   &\qquad+  \frac{1}{\Ca} \, \dn_1 \left( \partial_{Sx} \test{u}_{1x} \partial_{Sx} n_x  + 
                                                     \partial_{Sy} \test{u}_{1x} \partial_{Sy} n_x  + 
                                                     \partial_{Sx} \test{u}_{1y} \partial_{Sx} n_y  + 
                                                     \partial_{Sy} \test{u}_{1y} \partial_{Sy} n_y  
           \right)  \nonumber\\ 
    &\qquad-  \frac{1}{\Ca} \, [  n_x (\partial_{Sx} \test{u}_{1x} \partial_{Sx} \dn_{1} \!+\! \partial_{Sy} \test{u}_{1x} \partial_{Sy} \dn_{1}) + n_y (\partial_{Sx} \test{u}_{1y} \partial_{Sx} \dn_{1} \!+\! \partial_{Sy} \test{u}_{1y} \partial_{Sy} \dn_{1})  ] \Big{\}}  }{\Sigma_{B0}}  \nonumber\\ 
    +& \intS{ \test{u}_{1x} \Delta p_1   }{\Sigma_{in}} \,,
\end{align}
where variables $\tau_{1nx}$ and $\tau_{1ny}$ represent the stresses exerted on the wall in the $x$ and $y$ directions, respectively. The functions $u_{1x}$ and $u_{1y}$ are periodic, thus fulfilling the periodicity equations \eqref{perv1} and \eqref{perdnv1}. Multipliying the impermeability condition \eqref{imp1} by the test function of the degree of freedom for DBP $\test{\dn}_1$ and integrating over the boundary $\Sigma_{B0}$ leads to 
\begin{align}
0=& \intS{ \Big{[} \test{\dn}_1 (u_{1x} n_x \!+\! u_{1y} n_y) + ( u_{0x} \partial_{Sx} \test{\dn}_{1}  + u_{0y} \partial_{Sy} \test{\dn}_{1} )\dn_1 \Big{]} }{\Sigma_{B0}}   \,.    
\end{align}
In addition, the global equations \eqref{Flow1} and \eqref{Surface1} close the system of equations.

\subsection{Discussion}\label{Disc}

\begin{figure}[h]
 \centering
\begin{tikzpicture}[baseline,remember picture]
\begin{axis}[%
width=.62\textwidth,
height=0.1\textwidth,
xmin=0,
xmax=.62,
ymin=0,
ymax=1,
    colormap/jet,
    colorbar,
    point meta min=0,
    point meta max=.12,
    colorbar style={
        width=.2cm,
        ytick={0,.06,.12},
        yticklabels={$0.00$,$0.06$,$0.12$},
        },
        hide axis,
]
\addplot[thick,blue] graphics[xmin=0.0 ,ymin=0,xmax=0.3,ymax=1] {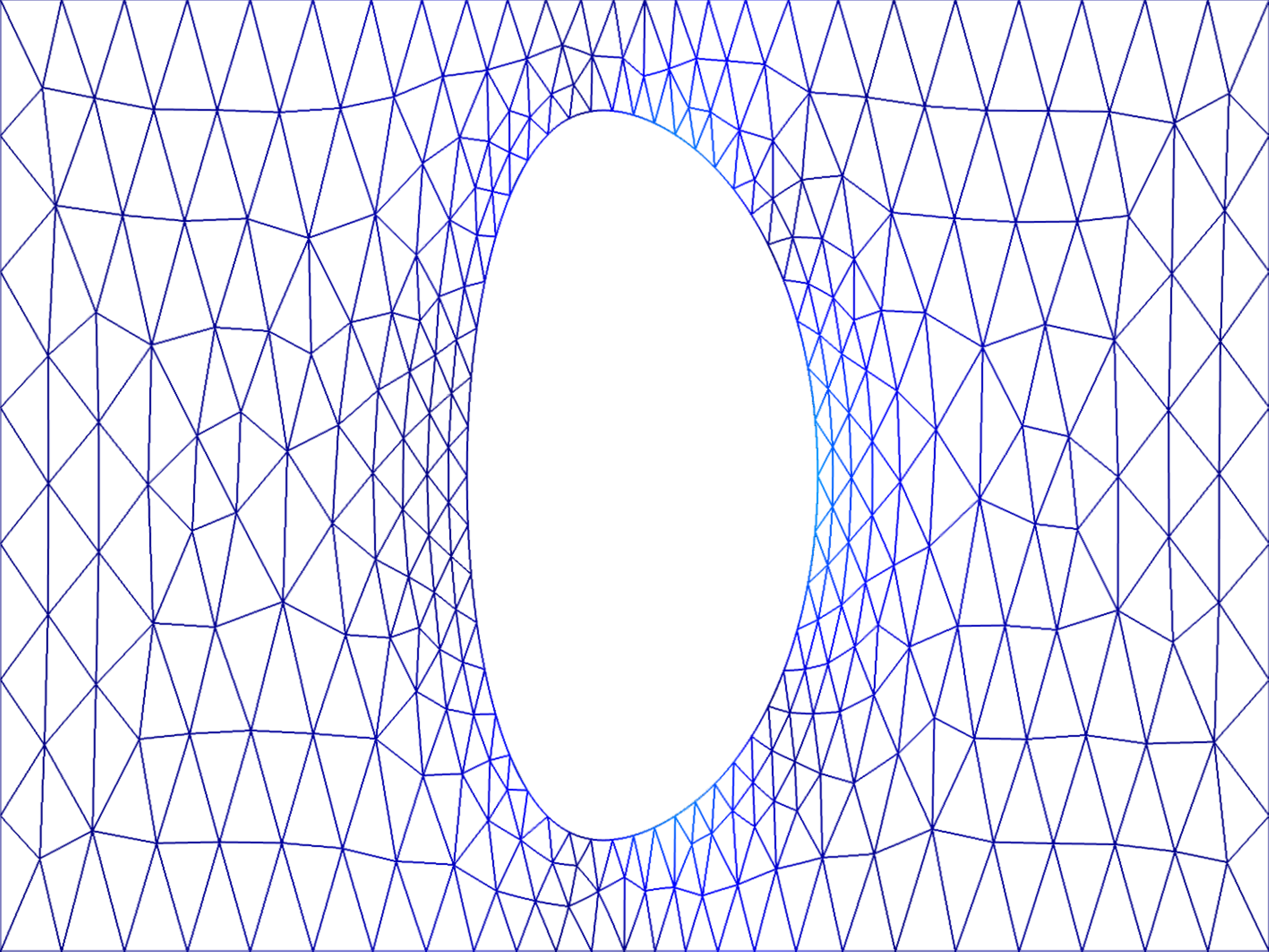};
\addplot[thick,blue] graphics[xmin=0.32,ymin=0,xmax=0.62,ymax=1] {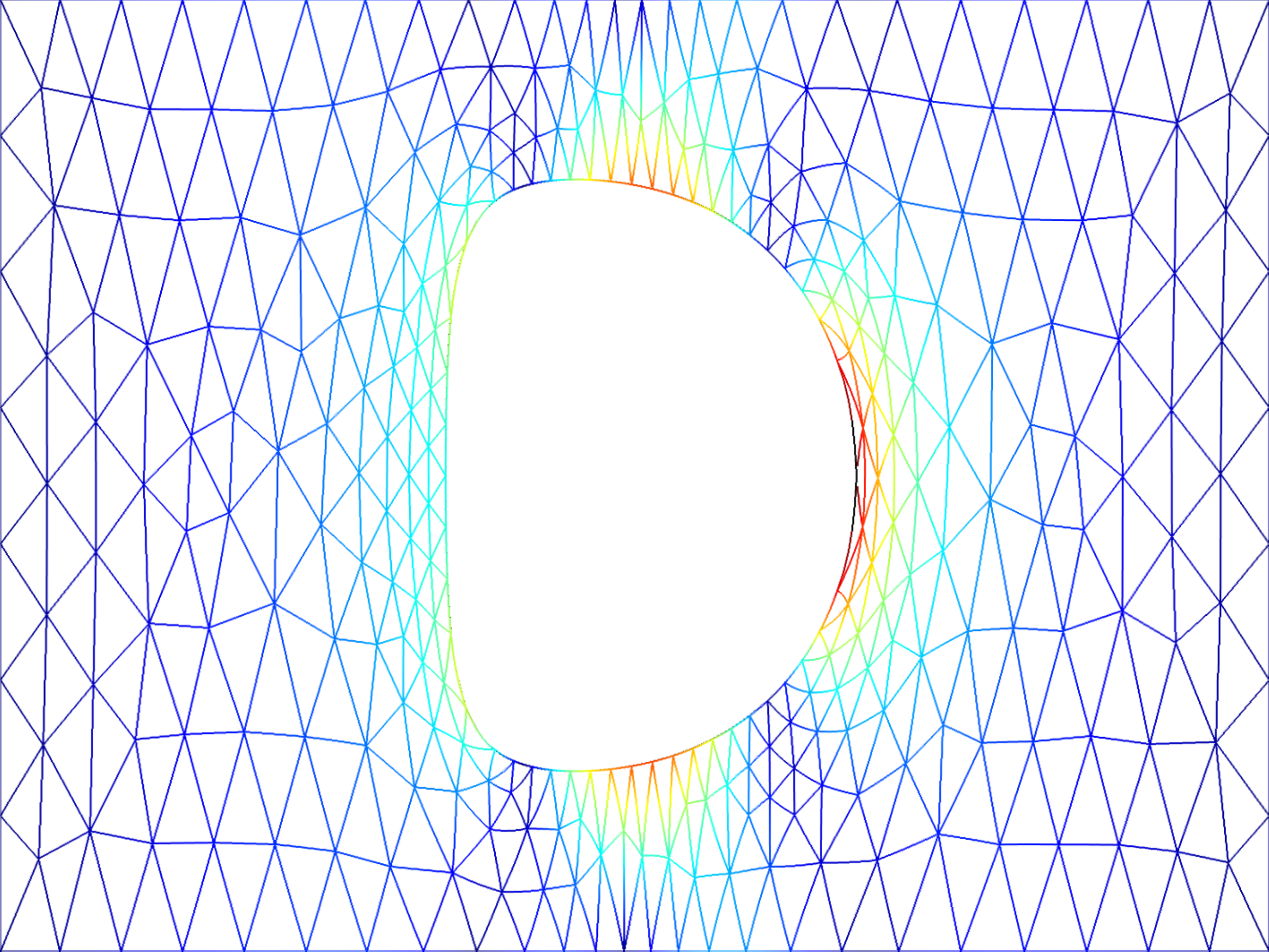};
%\node[anchor=south west] at (0,0) {\includegraphics[width=.3\textwidth,height=.1\textwidth]{./figures/Malla_R04_Ca01_v3.png}};
%\includegraphics[width=.3\textwidth,height=.1\textwidth]{./figures/Malla_R04_Ca05_v3.png}
\end{axis}
\end{tikzpicture}
\caption{Wireframe of the deformed mesh for a bubble of size $\V_B=\pi 0.4^2$ and $L=3$, for (left) $Ca=0.1$ and (right) $Ca=0.5$. Colormaps $|| \vx - \vX || $.}
\label{Wireframe}
\end{figure}

In \figref{Wireframe}, we depict the mesh for a centred bubble for two values of $\Ca$. We can observe how the mesh on the bubble interface remains equally spaced since the deformation is governed by Laplace equation within the boundary since the nonhomogeneous term in its governing equation \eqref{BALEeq} has no contribution within the boundary. This feature can be used to avoid remeshing in many situations.

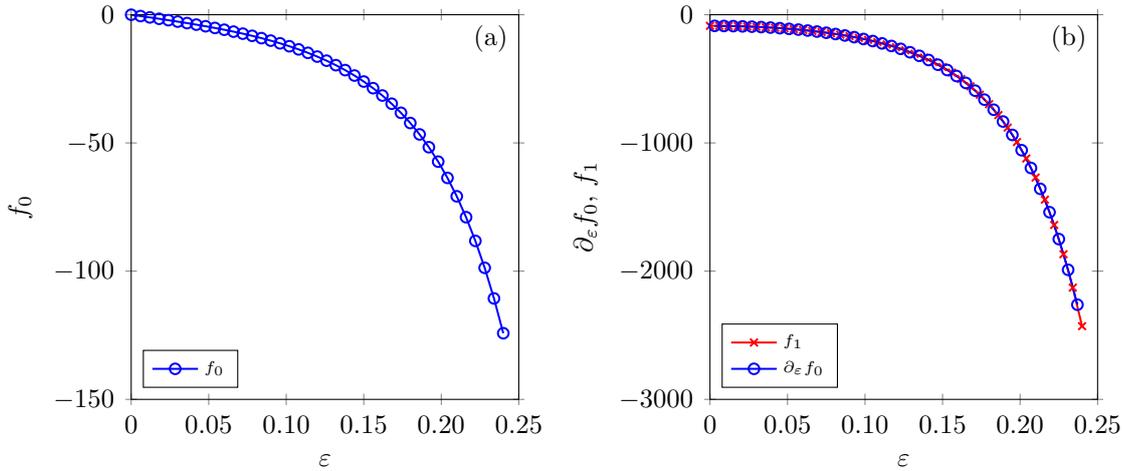
\begin{figure}[h]
% This file was created by matlab2tikz.
%
%The latest updates can be retrieved from
%  http://www.mathworks.com/matlabcentral/fileexchange/22022-matlab2tikz-matlab2tikz
%where you can also make suggestions and rate matlab2tikz.
%
\definecolor{mycolor1}{rgb}{0.00000,0.50000,1.00000}%
\definecolor{mycolor2}{rgb}{0.00000,1.00000,1.00000}%
\definecolor{mycolor3}{rgb}{1.00000,1.00000,0.00000}%
\begin{tikzpicture}[baseline,remember picture]

\begin{axis}[%
width=0.3\textwidth,
height=0.3\textwidth,
at={(0\textwidth,0\textwidth)},
scale only axis,
%
%xmode=log,
xmin=0,
xmax=.25,
xtick={0,0.05,0.1,0.15,.20,.25},
xticklabels={$0$,$0.05$,$0.10$,$0.15$,$0.20$,$0.25$},
%xminorticks=true,
xlabel={$\varepsilon$},
%
%ymode=log,
ymin=-150,
ymax=0,
ytick={-150,-100,-50,0},
yticklabels={$-150$, $-100$, $-50$, $0$},
%%%yminorticks=true,
ylabel={$f_0 $},
%
%axis x line*=bottom,
%axis y line*=left,
%legend style={legend cell align=left,row sep =-3pt,at={(1,0.025)},anchor=south east},
%legend pos=north west,
%legend style={at={(1.1,.5)},anchor=west}
legend pos=south west,
]

\node[anchor=north east] at (rel axis cs:1,1) {(a)};

\addplot [color=blue,mark=o, mark size=2pt]
  table[row sep=crcr]{%
     0.00000    -0.00192\\
     0.00600    -0.52425\\
     0.01200    -1.05120\\
     0.01800    -1.58490\\
     0.02400    -2.12980\\
     0.03000    -2.68890\\
     0.03600    -3.26760\\
     0.04200    -3.86880\\
     0.04800    -4.49820\\
     0.05400    -5.15980\\
     0.06000    -5.85620\\
     0.06600    -6.59640\\
     0.07200    -7.38440\\
     0.07800    -8.22720\\
     0.08400    -9.13060\\
     0.09000   -10.10200\\
     0.09600   -11.15100\\
     0.10200   -12.28700\\
     0.10800   -13.51900\\
     0.11400   -14.86100\\
     0.12000   -16.32300\\
     0.12600   -17.92000\\
     0.13200   -19.67600\\
     0.13800   -21.59900\\
     0.14400   -23.71700\\
     0.15000   -26.05700\\
     0.15600   -28.64400\\
     0.16200   -31.51200\\
     0.16800   -34.70400\\
     0.17400   -38.25800\\
     0.18000   -42.23200\\
     0.18600   -46.67800\\
     0.19200   -51.67600\\
     0.19800   -57.30100\\
     0.20400   -63.64100\\
     0.21000   -70.81300\\
     0.21600   -78.96100\\
     0.22200   -88.20300\\
     0.22800   -98.70300\\
     0.23400  -110.64000\\
     0.24000  -124.21000\\
};
\addlegendentry{$ f_0$};

\end{axis}
\end{tikzpicture}%
% This file was created by matlab2tikz.
%
%The latest updates can be retrieved from
%  http://www.mathworks.com/matlabcentral/fileexchange/22022-matlab2tikz-matlab2tikz
%where you can also make suggestions and rate matlab2tikz.
%
\definecolor{mycolor1}{rgb}{0.00000,0.50000,1.00000}%
\definecolor{mycolor2}{rgb}{0.00000,1.00000,1.00000}%
\definecolor{mycolor3}{rgb}{1.00000,1.00000,0.00000}%
\begin{tikzpicture}[baseline,remember picture]

\begin{axis}[%
width=0.3\textwidth,
height=0.3\textwidth,
at={(0\textwidth,0\textwidth)},
scale only axis,
%
%xmode=log,
xmin=0,
xmax=.25,
xtick={0,0.05,0.1,0.15,.20,.25},
xticklabels={$0$,$0.05$,$0.10$,$0.15$,$0.20$,$0.25$},
%xminorticks=true,
xlabel={$\varepsilon$},
%
%ymode=log,
ymin=-3000,
ymax=0,
ytick={-3000,-2000,-1000,0},
yticklabels={$-3000$, $-2000$, $-1000$, $0$},
%%%yminorticks=true,
ylabel={$\partial_\varepsilon f_0$, $f_1$},
%
%axis x line*=bottom,
%axis y line*=left,
%legend style={legend cell align=left,row sep =-3pt,at={(1,0.025)},anchor=south east},
legend pos=south west,
]

\node[anchor=north east] at (rel axis cs:1,1) {(b)};

\addplot [color=red,mark=x, mark size=2pt]
  table[row sep=crcr]{%
  0.0000e+00  -8.6748e+01\\
   6.0000e-03  -8.7031e+01\\
   1.2000e-02  -8.7897e+01\\
   1.8000e-02  -8.9404e+01\\
   2.4000e-02  -9.1550e+01\\
   3.0000e-02  -9.4290e+01\\
   3.6000e-02  -9.7751e+01\\
   4.2000e-02  -1.0187e+02\\
   4.8000e-02  -1.0677e+02\\
   5.4000e-02  -1.1245e+02\\
   6.0000e-02  -1.1905e+02\\
   6.6000e-02  -1.2649e+02\\
   7.2000e-02  -1.3490e+02\\
   7.8000e-02  -1.4450e+02\\
   8.4000e-02  -1.5522e+02\\
   9.0000e-02  -1.6726e+02\\
   9.6000e-02  -1.8072e+02\\
   1.0200e-01  -1.9595e+02\\
   1.0800e-01  -2.1296e+02\\
   1.1400e-01  -2.3190e+02\\
   1.2000e-01  -2.5331e+02\\
   1.2600e-01  -2.7727e+02\\
   1.3200e-01  -3.0446e+02\\
   1.3800e-01  -3.3468e+02\\
   1.4400e-01  -3.6881e+02\\
   1.5000e-01  -4.0788e+02\\
   1.5600e-01  -4.5120e+02\\
   1.6200e-01  -5.0131e+02\\
   1.6800e-01  -5.5843e+02\\
   1.7400e-01  -6.2343e+02\\
   1.8000e-01  -6.9756e+02\\
   1.8600e-01  -7.8273e+02\\
   1.9200e-01  -8.8026e+02\\
   1.9800e-01  -9.9292e+02\\
   2.0400e-01  -1.1214e+03\\
   2.1000e-01  -1.2700e+03\\
   2.1600e-01  -1.4425e+03\\
   2.2200e-01  -1.6398e+03\\
   2.2800e-01  -1.8682e+03\\
   2.3400e-01  -2.1284e+03\\
   2.4000e-01  -2.4289e+03\\
};
\addlegendentry{$f_1$};

\addplot [color=blue,mark=o, mark size=2pt]
  table[row sep=crcr]{%
   3.0000e-03  -8.7054e+01\\
   9.0000e-03  -8.7825e+01\\
   1.5000e-02  -8.8950e+01\\
   2.1000e-02  -9.0817e+01\\
   2.7000e-02  -9.3183e+01\\
   3.3000e-02  -9.6450e+01\\
   3.9000e-02  -1.0020e+02\\
   4.5000e-02  -1.0490e+02\\
   5.1000e-02  -1.1027e+02\\
   5.7000e-02  -1.1607e+02\\
   6.3000e-02  -1.2337e+02\\
   6.9000e-02  -1.3133e+02\\
   7.5000e-02  -1.4047e+02\\
   8.1000e-02  -1.5057e+02\\
   8.7000e-02  -1.6190e+02\\
   9.3000e-02  -1.7483e+02\\
   9.9000e-02  -1.8933e+02\\
   1.0500e-01  -2.0533e+02\\
   1.1100e-01  -2.2367e+02\\
   1.1700e-01  -2.4367e+02\\
   1.2300e-01  -2.6617e+02\\
   1.2900e-01  -2.9267e+02\\
   1.3500e-01  -3.2050e+02\\
   1.4100e-01  -3.5300e+02\\
   1.4700e-01  -3.9000e+02\\
   1.5300e-01  -4.3117e+02\\
   1.5900e-01  -4.7800e+02\\
   1.6500e-01  -5.3200e+02\\
   1.7100e-01  -5.9233e+02\\
   1.7700e-01  -6.6233e+02\\
   1.8300e-01  -7.4100e+02\\
   1.8900e-01  -8.3300e+02\\
   1.9500e-01  -9.3750e+02\\
   2.0100e-01  -1.0567e+03\\
   2.0700e-01  -1.1953e+03\\
   2.1300e-01  -1.3580e+03\\
   2.1900e-01  -1.5403e+03\\
   2.2500e-01  -1.7500e+03\\
   2.3100e-01  -1.9895e+03\\
   2.3700e-01  -2.2617e+03\\
};
\addlegendentry{$\partial_\varepsilon f_0$};

\end{axis}
\end{tikzpicture}%
\caption{Validation of the procedure for a bubble of size $\V_B=\pi 0.2^2$, $Ca=0.2$ and $L=3$. Value of (a) the migration force $f_0$ and (b) its sensitivity to a variation in $\varepsilon$ obtained as either its derivative with respect to the position $\partial_\varepsilon f_0$ or its perturbation $f_1$.}
\label{Validation}
\end{figure}

The sensitivity of the migration force $f$ to changes in the equilibrium position $\varepsilon$ can be obtained by taking the derivative for the zeroth-order solution, based on BALE method, or from the first-order solution that results from the perturbation of the equilibrium position $\varepsilon$, based on the DBP method, as schematised in \figref{Sketch1}. The agreement between the results obtained from both methods, BALE and DBP, validates them. For this purpose, we consider that the perturbation of the function $f(\varepsilon)$ writes, due to the asymptotic expansion \eqref{expansb}, as
\begin{align}
f (\varepsilon + \epsilon) \approx f_0 (\varepsilon) + \epsilon f_1(\varepsilon) \,,
\end{align}
and the function $f_0 (\varepsilon)$ writes, after Taylor expansion, as
\begin{align}
f_0 (\varepsilon + \epsilon) \approx f_0 (\varepsilon) + \epsilon \frac { \partial f_0}{\partial \varepsilon}\Big{\vert}_\varepsilon \,.
\end{align}
Then, since $f (\varepsilon + \epsilon) = f_0 (\varepsilon + \epsilon)$, $f_1$ is the sensitivity of the migration force to changes in the equilibrium position
\begin{align}
\label{Sensit}
f_1 =  \frac { \partial f_0}{\partial \varepsilon}\Big{\vert}_\varepsilon \,.
\end{align}
In \figref{Validation}, we depict the functions $f_0(\varepsilon)$ and $f_1(\varepsilon)$. We compute $\partial_{\varepsilon} f_0$ applying centred finite difference to the discrete data for $f_0(\varepsilon)$ depicted in \figref{Validation}a. In \figref{Validation}b, we depict the computed $\partial_{\varepsilon} f_0$ and $f_1$ and we can observe that \eqref{Sensit} perfectly holds, thus validating both methods.  

It is worth noting that, to obtain the sensitivity of the migration force against the position, the linearisation is accessory but very useful for the presentation of the method and validation. However, its use is very convenient in other problems such as in \cite{Rivero2018a} where it reduces the computational cost with respect to the nonlinear deformable domain since the perturbation system of equations is linear in the former case. Furthermore, in some situations such as stability analyses or steady problems, this procedure makes it possible and avoids the need of full transient analysis. Concerning shape optimization, the optimal modification of the shape reduces the computational cost of remeshing to obtain the optimal deformation of the geometry.

 %%%%%%%%%%%%%%%%%%%%%%%%%%
%%%%%%%%%%%%%%%%%%%%%%%%%%
\section{Conclusions}\label{Conclusions}

In this paper, we propose two methods to treat partial differential equations (PDE) defined on deformable domains for geometries that undergo either large or small deformations. On the one hand, for large deformations of the domain, we propose the Boundary Arbitrary Lagrangian-Eulerian (BALE) method to track the boundary of the domain analogue and complementary to the Arbitrary Lagrangian-Eulerian method (ALE). It is analogue in the sense that both the BALE and ALE methods rely on a change of variable which fulfils a given PDE. However, the BALE method relies on only one degree of freedom for the displacement of the boundary. It is complementary in the sense that the BALE method can be used for the displacement of the boundary, being the boundary condition of the PDE used in the change of variable for the ALE method. On the other hand, for small deformations of the domain, the domain can be perturbed and the equations and boundary conditions can be written in the unperturbed domain and boundary which are known. We named this method the deformable boundary perturbation (DBP) method. These two methods complement the existing tools available up to date and fill the needs of systematic tool to treat free interfaces in a non case-dependent manner.
 
The BALE method has several advantages that should be noted. First, it is as systematic and not case-dependent as the ALE method and therefore, it can be applied for any geometrical configuration with structured and unstructured meshes. Second, it reduces the distortion of the mesh at the boundaries. As we observe in the example of a bubble in a microchannel flow, an equally spaced mesh at the boundary remains equally spaced after deformation. And third, it removes the two degrees of freedom of the mesh within the boundary, which yields to no variation of the boundary, whereas the degree of freedom out of the boundary remains free. It is very useful for stationary analysis in the presence of deformable domains where the Eulerian description is hardly avoidable.

The DBP method is a tool to treat unknown boundaries that makes possible or facilitates stability analysis and expansions in terms of a small parameters as well as shape optimisation with gradient-based methods. It is a systematic and avoid case-dependent methods which may discourage the application of stability analysis and expansions in terms of small parameters due to non-simple geometries. Concerning the shape optimisation, different shape modifications can be tested in a very accurate manner by solving a linear problem obtained from the exact linearisation possible by the use of DBP. It decreases the computational cost with respect to shape modification that leads to remeshing the modified geometry while preserving the accuracy. 

Despite the uncommon use of the boundary exterior differential operator, we found it very useful for several reasons. First, in the case of DBP applied to mixed boundary conditions, the flux through the generatrix of a quantity can be expressed in terms of this operator thanks to the boundary Stokes theorem. Second, surface tension is expressed in terms of this operator in a very compact manner, and its perturbation can also be written in terms of this operator. And third, it can be very easily implemented using finite element methods, but it is not restricted to it and its definition in Cartesian coordinates can be used instead.
 
We have applied both methods to a fluid mechanical system with free interface and detailed the implementation in weak form. We have then validated the results by comparison between the results obtained from the two methods. The presented methods can be applied to other physical systems with deformable domains and governed by PDEs, provided that the boundaries of the domain are sharp and therefore, not compatible with diffuse interface models.

%%%%%%%%%%%%%%%%%%%%%%%%%%
%%%%%%%%%%%%%%%%%%%%%%%%%%
\section*{Acknowledgements}

We thank the Brussels region for the financial support of this project through the WBGreen-MicroEco project. We also thank the F.R.S.-FNRS for financial support through the WOLFLOW project as well as the IAP-7/38 MicroMAST project for supporting this research.

\appendix

\section{Perturbation of $\hgradS$}\label{AppA}

In this appendix, we provide an alternative procedure to obtain the perturbation of $\delta (\vnS \dd \Gamma)$. For this purpose, let us consider the unitary vector tangent to the curve $\Gamma$, $ \vt = \vn \x \vnS $ as well as the line vector $\vt \dd \Gamma$, which according to \cite{batchelor2000introduction} can be perturbed as
\begin{align}
\label{pertdr}
\delta (\vt \dd \Gamma) =  \vt \dd \Gamma \cdot \gradS \vdn \,. 
\end{align}
The perturbation of $\vnS \, \dd \Gamma = \vt \dd \Gamma \x \vn $ writes, using the product rule and the equations \eqref{pertdr} and \eqref{deltan}, as 
\begin{align}
\label{eqap1}
\delta ( \vnS \, \dd \Gamma) =  \delta (\vt \dd \Gamma)  \x \vn + \vt \dd \Gamma \x  \delta (\vn) =  \vt \dd \Gamma \cdot \gradS \vdn  \x \vn - (\vn \x \vnS) \dd \Gamma \x \gradS \vdn \cdot \vn 
\end{align}
The first term of \eqref{eqap1} can be simplified by projecting with the identity tensor on the right. After the circular shift property of the triple product, it can be written as
\begin{align}
 \vt \dd \Gamma  \cdot \gradS \vdn  \x \vn =  \dd \Gamma \vt \cdot \gradS \vdn  \cdot ( \vn  \x \id)  
\end{align}
Since, $\id = \vt\vt + \vnS\vnS + \vn\vn$, the vectorial product can be written as $\vn \x \idS =  \vt \vnS - \vnS \vt $. Adding and subtracting $\dd \Gamma \vnS \vnS \cdot \gradS \vdn  \cdot  \vnS $ leads, after convenient rearrangement of positive and negative terms, to 
\begin{align}
 \vt \dd \Gamma  \cdot \gradS \vdn  \x \vn =  \dd \Gamma \vnS (\vt \cdot \gradS  \vdn \cdot \vt + \vnS \cdot \gradS  \vdn \cdot \vnS ) 
    - \dd \Gamma  [ (\vt \vt+ \vnS \vnS) \cdot \gradS \vdn  \cdot  \vnS ] 
\end{align}
Furthermore, taking into account that $\vt \cdot \gradS  \vdn \cdot \vt + \vnS \cdot \gradS  \vdn \cdot \vnS  = \gradS \cdot \vdn$ and $(\vt \vt+ \vnS \vnS) \cdot \gradS = \gradS$, the first term of the perturbation of $\vnS \, \dd \Gamma $ can be written as
\begin{align}
\label{first}
 \vt \dd \Gamma                \cdot \gradS \vdn  \x \vn =  \dd \Gamma \vnS \cdot [ \idS \gradS \cdot \vdn   - (\gradS \vdn)^T  ]    
\end{align}

The second term of the RHS of \eqref{eqap1} can be rewritten using the circular shift property of the triple vector product
\begin{align}
\label{second}
-(\vn \x \vnS) \dd \Gamma \x \gradS \vdn \cdot \vn = \vnS \dd \Gamma \cdot (\gradS \vdn) \cdot \vn \vn \,.
\end{align}

The substitution of \eqref{first}-\eqref{second} leads to the same equation as the geometrical method \eqref{eq32}.

\bibliographystyle{plain}
\bibliography{mybib.bib}
%\begin{thebibliography}{00}
%
%%% \bibitem[Author(year)]{label}
%%% Text of bibliographic item
%
%\bibitem[ ()]{}
%
%\end{thebibliography}
\end{document}